\documentclass[iicol]{sn-jnl}% Basic Springer Nature Reference Style/Chemistry Reference Style

\usepackage{graphicx}%
\usepackage{multirow}%
\usepackage{multicol}
\usepackage{amsmath,amssymb,amsfonts}%
\usepackage{amsthm}%
\usepackage{mathrsfs}%
\usepackage[title]{appendix}%
\usepackage{xcolor}%
\usepackage{textcomp}%
\usepackage{manyfoot}%
\usepackage{booktabs}%
\usepackage{algorithm}%
\usepackage{algorithmicx}%
\usepackage{algpseudocode}%
\usepackage{listings}%

\usepackage{cite}
\usepackage{lipsum}   %用于占位
\usepackage{array}

\usepackage{makecell}
\usepackage{multirow}
\usepackage[square,numbers]{natbib} % 数字化citation
\usepackage{tabularx}
\usepackage{adjustbox}

\begin{document}

\title[Article Title]{Mixed Precision Block-Jacobi Preconditioner: Algorithms, Performance Evaluation and Feature Analysis}

\author[1,2]{\fnm{Ningxi} \sur{Tian}}\email{tianningxi22@gscaep.ac.cn}

\author[2]{\fnm{Silu} \sur{Huang}}\email{huang\_silu@iapcm.ac.cn}

\author*[2]{\fnm{Xiaowen} \sur{Xu}}\email{xwxu@iapcm.ac.cn} 

\affil[1]{\orgname{Graduate School of Chinese Academy of Engineering Physics}, \orgaddress{\city{Beijing}, \postcode{100094}, \country{China}}}

\affil[2]{\orgdiv{Laboratory of Computational Physics}, \orgname{Institute of Applied Physics and Computational Mathematics}, \orgaddress{\city{Beijing}, \postcode{100088}, \country{China}}} 

\abstract{
In this paper, we propose two mixed precision algorithms for block-Jacobi preconditioner(BJAC): a fixed low precision strategy and an adaptive precision strategy. We evaluate the performance improvement of the proposed mixed precision BJAC preconditioners combined with the preconditioned conjugate gradient(PCG) method using problems including diffusion equations and radiation hydrodynamics equations. Numerical results show that, compared to the uniform high precision PCG, the mixed precision preconditioners can achieve speedups from $1.3\times$ to $1.8\times$ without lossing accuracy. Furthermore, we observe the phenomenon of convergence delay in some test cases for the mixed precision preconditioners, and analyse the correlation between matrix features and convergence delay behaviors, some interesting conclusions are obtained which are significant and valuable for the design of more efficient mixed precision preconditioners.
}


\maketitle

%% ---------------------------------------------------------------------------
\section{Introduction}\label{sec:1}

Solving large-scale sparse linear systems arising from the discretization of partial differential equations is a crucial challenge in scientific and engineering computing. In many realistic simulations such as laser fusion, complex fluid dynamics, structural mechanics, and electronic systems, the solution of such systems has become a main performance bottleneck \cite{Xu_PAMG}. 

The preconditioned Krylov subspace methods are widely used for solving sparse linear systems in practical applications \cite{Saad}. The main idea of this method is to transform the original system $Ax=b$ into an equivalent preconditioned systems $M^{-1}Ax = M^{-1}b$ via constructing a preconditioner $M^{-1}$, which is an approximation of $A^{-1}$, then the preconditioned systems is solved using Krylov subspace methods. The aim of preconditioiner is to accelerate the convergence speed of Krylov subspace iterations. The typical preconditioners include the (block) Jacobi method \cite{Saad}, incomplete LU decomposition method \cite{Saad,Golub}, multigrid \cite{Saad,Briggs_MG,Xu_AMG} and domain decomposition methods \cite{Saad,Toselli_DDM}, etc. These preconditioners are widely used in  practical applications and usually dominate the computation time, thus improving their efficiency is an important issue. 
In practical applications, the efficiency of preconditioners are both problem- and architecture-specific. 
On one hand, it is crucial to design the problem-specific algorithms to improve the convergence rate. On the other hand, performance optimization for target architecture is required to improve the floating-point operation efficiency(Flops).

In recent years, mixed precision has become a new technique to improve the performance of algorithms. 
Specifically, mixed precision technique employs low precision in some components of the algorithm without losing the accuracy required by the application, which has less computation cost, memory storage, data movement overheads and lower energy consumption. 
Currently, the mixed precision preconditioner has been investigated for various numerical algorithms \cite{Higham_MP,Abdelfattah_MP,Dongarra_MP}. 

The preconditioner, as mentioned above, is essentially an approximation to the inverse of the original matrix. Thus, a straightforward idea is designing a mixed precision preconditioner using low precision. The primary advantage of low precision preconditioners is the computational cost reduction of single iteration. However, it may lead to a decrease in convergence speed or even failure in convergence. The key to achieving performance improvement is the trade-off between the computational efficiency of single iteration and the speed of convergence, which needs to be designed carefully based on the application features and the algorithm characteristics. 

In this paper, we consider the block-Jacobi(BJAC) preconditioner, which is a kind of basic preconditioiner suitable for massively parallel computations \cite{Saad, Dubois_JAC}. We propose two mixed precision strategies for BJAC preconditioners, fixed-mixed-precision based BJAC (fMP-BJAC) and adaptive-mixed-precision based BJAC (aMP-BJAC).
When applied to the Krylov subspace iterative methods, the key is using the low precision format in BJAC and employing the high precision format in other operations in Krylov iterations. 
During the iterative process, the fMP-BJAC preconditioner uses a fixed low precision for both storage and computation, while the aMP-BJAC preconditioner tune its precision automatically based on the convergence behavior. 
Experimental results of the three-dimensional diffusion equations and radiation hydrodynamics equations, show the performance improvements without sacrificing accuracy. 
Moreover, we also investigate that the proposed mixed precision preconditioners may lead to the convergence delays of the Krylov subspace method, which associate with the application features, such as the multiscale and diagonal dominance. %[Haoyue]

The rest of the paper is organized as follows. A brief overview on current works of mixed precision preconditioners is introduced in Section \ref{sec:2}. Section \ref{sec:3} introduces the block-Jacobi preconditioner and mixed precision preconditioners we proposed. Numerical results and performance evaluation as well as feature analysis for the proposed preconditioners are presented in Section \ref{sec:4}. Finally, we summarise the conclusions and the issues for further research in Section \ref{sec:5}. 
\section{Related Work}\label{sec:2}

For mixed precision preconditioners, a main strategy is to integrate low precision preconditioners into high precision (usually double precision) Krylov subspace methods  \cite{Giraud_MPDDM, Arioli_MPfgmres, Kronbichler_MP, Anzt_aMPJac, Flegar_aMPJac, Carson_MPSPI,Chu_MPSPI, Zhang_MPSPI}. Giraud et al.(2008) \cite{Giraud_MPDDM} proposed a single precision domain decomposition preconditioner combined with a double precision CG algorithm, the results on 2D and 3D elliptic equations illustrate the improvement in reducing CPU time and storage as well as communication costs with iteration counts increasing slightly. 
Arioli et al.(2009) \cite{Arioli_MPfgmres} theoretically and experimentally confirmed that a single precision LU preconditioner combined with a double precision FGMRES(Flexible Generalized Minimal Residual) method maintains backward stability akin to double precision. 
Kronbichler et al.(2019) \cite{Kronbichler_MP} investigated a single precision geometric multigrid preconditioner with a double precision CG algorithm on GPUs, demonstrating comparable discretization error to double precision while accelerating computations by 47\% to 83\% for the Laplace problem. 
Anzt et al.(2019) \cite{Anzt_aMPJac} introduced an adaptive block-Jacobi preconditioner utilizing mixed precision storage. This approach adjusts the storage precision of each diagonal block's inverse matrix based on its condition number among double precision, single precision, half precision, and other customized low precision formats, while maintaining double precision in results. The numerical results for the SuiteSparse matrix set demonstrate that this adaptive mixed precision storage preconditioner reduces PCG algorithm solution times by 10\% to 30\% compared to algorithms using double precision exclusively \cite{Flegar_aMPJac}.
Carson and Khan (2023) \cite{Carson_MPSPI} analysed the sparse approximate inverses in finite precision and investigated the behaviour of the sparse approximate inverse preconditioner within five-precision GMRES-IR refinement.
Chu (2023) \cite{Chu_MPSPI} designed a mixed precision static sparse approximation inverse preconditioner and a mixed precision heuristic sparse approximation inverse preconditioner on GPU. 
Zhang et al.(2024) \cite{Zhang_MPSPI} propose and investigate a single precision block sparse approximation inverse preconditioner on Tensor core.

The current works on mixed precision preconditioners show the potential in improving the performance for solving sparse linear systems. However, on one hand, the test problems in current work mainly arising from the model problems or the problems from SuiteSparse matrix set with small scale, on the other hand, the impact on convergence behavior of the mixed preconditioners, which is an important aspect that impacts on overall performance, is rarely be concerned in current work. 
Further research is required on these issues.
\section{Mixed Precision Block-Jacobi Preconditioner}\label{sec:3}

\subsection{The Block-Jacobi Preconditioner(BJAC)}

For a sparse matrix $A$, considering its block partitioning  $A={(A_{ij})}_{nb\times nb} \in R^{n\times n}$, where $nb$ is the number of blocks, $A_{ij} \in R^{n_s \times n_s}$ is the sub-matrix of size $n_s=n/nb$. Let $D_b=diag\{A_{ii},i=1,2,\dots,nb\}$ denotes the diagonal block matrix, then the block-Jacobi preconditioner(BJAC) is defined as follows \cite{Dubois_JAC},

\begin{equation}
	\label{eq3-1}
	M^{-1}_{bjac} = \sum^{k-1}_{j=0}{(I-D^{-1}_b A)}^j D^{-1}_b
\end{equation}
where $I$ is the identity matrix, $k$ is the number of iterations. In parallel computing, the number of blocks $nb$ is usually set to  the number of processes.

Since the cost of the exact inverse matrices $D^{-1}_b$ is expensive due to the high complexity of $A^{-1}_{ii}$, it is usually to use an approximation version  $\hat{D}^{-1}_b$ instead, then an approximation of the block-Jacobi preconditioner is obtained, which as defined as follows,

\begin{equation}
	\label{eq3-2}
	\hat{M}^{-1}_{bjac} = \sum^{k-1}_{j=0}{(I-\hat{D}^{-1}_bA)}^j\hat{D}^{-1}_b
\end{equation}
In particular, if the approximate inverse $\hat{A}^{-1}_{ii}$ of $A^{-1}_{ii}$ is implemented via a Jacobi iterative method $\hat{A}^{-1}_{ii} = \sum^{t-1}_{i=0}{(I-{D}^{-1}_{ii} A_{ii})}^iD^{-1}_{ii}$, where $D_{ii}$ is diagonal matrix of $A_{ii}$, $t$ is the number of iterations within blocks, then $\hat{D}^{-1}_b$ has the following form: 

\begin{equation}
	\label{eq3-3}
	\resizebox{0.9\hsize}{!}{$
	\hat{D}^{-1}_b = 
	\left[
	\begin{array}{ccc}
		\displaystyle\sum^{t-1}_{i=0}{(I-{D}^{-1}_{11} A_{11})}^iD^{-1}_{11} & & \\
		& \ddots  & \\
		&  &\displaystyle\sum^{t-1}_{i=0}{(I-{D}^{-1}_{nb,nb}A_{nb,nb})}^i D^{-1}_{nb,nb}
	\end{array}
	\right]
$}
\end{equation}

Although the $\hat{M}^{-1}_{bjac}$ is an approximation of $M^{-1}_{bjac}$ , and the latter can be regarded as a special case of the former, 
for convenience,  we still calling the approximation version in (\ref{eq3-2}) the block-Jacobi preconditioner(BJAC) in this paper, and still using $M^{-1}_{bjac}$ instead of $\hat{M}^{-1}_{bjac}$ to denote approximated BJAC preconditioner in the following sections. 

As shown in (\ref{eq3-2}) and (\ref{eq3-3}), the BJAC preconditioner has two parameters: the number of inter-block iterations $k$ and the number of intra-block iterations $t$. Then, by choosing $k$ and $t$, the preconditioners have different approximations to $A^{-1}$. In general, the larger $t$ and $k$ are, the better the BJAC approximation to $A^{-1}$ is, 
hence the better it accelerates the convergence of the Krylov subspace methods, meanwhile, the cost of the preconditioner is more expensive. Therefore, these two parameters need to be chosen carefully in practical problems in order to obtain the optimal performance. In particular, the BJAC preconditioner degenerates to the diagonal scaling preconditioner $D^{-1}$ if $k=t=1$. 

\subsection{Mixed Precision Block-Jacobi Preconditioner(MP-BJAC)}\label{sec:3.2}

Although the construction of BJAC is straightforward, it is one of the most classical and basic preconditioners, which is used directly or indirectly in many practical applications since it is suitable for massively parallel computation \cite{Saad}. Therefore, it is a typical  example used to verify the potential and efficiency of the mixed precision preconditioner. 

For the Krylov subspace methods, the behaviour of the preconditioner $M^{-1}$ is reflected in the operation of matrix-vector multiplication $z=M^{-1}r$ where $r$ is a vector produced during the iteration process of Krylov methods.
In this paper, we taking the Preconditioned Conjugate Gradient (PCG) method \cite{Hestenes_CG} as an example to design and verify the mixed precision BJAC preconditioner. Concretely, the low precision strategy is designed for the preconditioner, while the working precision is used for the remain operations in PCG method.

We define two precision formats: the high precision (i.e., working precision) format $p_w$ and the low precision format $p_l$. Denote the low precision BJAC preconditioners in which all floating-point computations are implemented in low precision format to be $M^{-1}_{bjac}(p_l)$. In this paper, we present two mixed precision BJAC preconditioner: the fixed low precision BJAC preconditioner(fMP-BJAC) and the adaptive precision BJAC preconditioner (aMP-BJAC).

\subsubsection{Fixed Low Precision BJAC Preconditioner}

For fMP-BJAC, as shown in Algorithm \ref{Algorithm 3.1}, before the preconditioner $M^{-1}_{bjac}(p_l)$ is multiplied with the vector $r$, the vector $r$ need to be converted from high to low precision (denoted as $\hat{r}$), and then compute the product of the preconditioner and the vector $\hat{r}$ in the low precision format, i.e., $\hat{z}=M^{-1} _{bjac}(p_l)\hat{r}$, and finally convert the low precision $\hat{z}$ to high precision $z$.

% \begin{algorithm}
% 	\caption{fMP-BJAC preconditioner }
%     % \setstretch{0.95}
% 	\SetAlgoLined
% 	\KwIn{$M^{-1}_{bjac}(p_l),r,p_w,p_l$}
% 	\KwOut{$z$}
% 	\label{Algorithm 3.1}
	
% 	$r = \text{ConvertPrecision}(\hat{r}, p_l)$ \\
% 	$\hat{z}=M^{-1}_{bjac}(p_l)\hat{r}$ \,  \textbackslash \textbackslash Compute in precision $p_l$ \\ 
% 	$z = \text{ConvertPrecision}(\hat{z}, p_w)$ \\
% 	return $z$\\
	
% \end{algorithm}

\begin{algorithm}
\caption{fMP-BJAC preconditioner}
\label{Algorithm 3.1}
\begin{algorithmic}[1] % [1] 表示行号
\Require $M^{-1}_{bjac}(p_l), r, p_w, p_l$
% \Ensure $z$
\State $r = \text{ConvertPrecision}(\hat{r}, p_l)$
\State $\hat{z} = M^{-1}_{bjac}(p_l) \cdot \hat{r}$ \Comment{Compute in precision $p_l$}
\State $z = \text{ConvertPrecision}(\hat{z}, p_w)$
\State \Return $z$
\end{algorithmic}
\end{algorithm}

The fMP-BJAC uses a low precision format for storing and computing, thereby reducing the time and data movement overhead for each iteration of the PCG algorithm. However, this mixed precision preconditioner may result in slower convergence of the PCG in some test problems, i.e., convergence delays occur, as discussed in the next section. Therefore, we further propose the adaptive precision BJAC (aMP-BJAC) preconditioner.

\subsubsection{Adaptive Precision BJAC Preconditioner}

For aMP-BJAC preconditioner, the precision is dynamically tuned between high and low precision based on the behavior of residual reduction during the iteration process. 

Two specific algorithms, aMP-BJAC(hl) and aMP-BJAC(lh), are defined respectively based on the order of switching ``high precision to low precision(hl)'' and ``low precision to high precision (lh)''. 
Algorithm \ref{Algorithm 3.2} and Algorithm \ref{Algorithm 3.3} illustrate the procedures of these two adaptive mixed precision preconditioners respectively. Taking the aMP-BJAC(hl) algorithm as an example, given the adaptive threshold $adp_{tol}$, the high precision is used when the norm of the relative residual ($relres$) is not less than $adp_{tol}$ (generally the early iteration stage), and the low precision is used when $relres$ is less than $adp_{tol}$ (generally the post-iteration stage).

% \begin{algorithm}
% 	\caption{aMP-BJAC(hl) preconditioner}
% 	\label{Algorithm 3.2}
% 	\SetAlgoLined
% 	\KwIn{$M^{-1}_{bjac},M^{-1}_{bjac}(p_l),r,p_w,p_l,relres,adp_{tol}$}
% 	\KwOut{$z$}
% 	\If{$relres \geq adp_{tol}$}{
% 		$z = M^{-1}_{bjac} r$
% 	}
% 	\Else{
% 		$z = \text{fMP-BJAC}(M^{-1}_{bjac}(p_l),r,p_w,p_l)$
% 	}
% 	\Return $z$
% \end{algorithm}

\begin{algorithm}
\caption{aMP-BJAC(hl) preconditioner}
\label{Algorithm 3.2}
\begin{algorithmic}[1] % [1] 表示行号
\Require $M^{-1}_{bjac}, M^{-1}_{bjac}(p_l), r, p_w, p_l, relres, adp_{tol}$
% \Ensure $z$
\If{$relres \geq adp_{tol}$}
    \State $z = M^{-1}_{bjac} \cdot r$
\Else
    \State $z = \text{fMP-BJAC}(M^{-1}_{bjac}(p_l), r, p_w, p_l)$
\EndIf
\State \Return $z$
\end{algorithmic}
\end{algorithm}

% \begin{algorithm}
% 	\caption{aMP-BJAC(hl) preconditioner}
% 	\SetAlgoLined
% 	\KwIn{$M^{-1}_{bjac},M^{-1}_{bjac}(p_l),r,p_w,p_l,relres,adp_{tol}$}
% 	\KwOut{$z$}
% 	\label{Algorithm 3.2}
	
% 	\eIf{ $relres \geq adp_{tol} $}{
% 		$ z = M^{-1}_{bjac }r$ 
% 	}{
% 		$z = \text{fMP-BJAC(}M^{-1}_{bjac}(p_l),r,p_w,p_l \text{)} $ 
% 	}
% 	return $z$
	
% \end{algorithm}

% \begin{algorithm}
% 	\caption{aMP-BJAC(lh) preconditioner}
	
% 	\SetAlgoLined
% 	\KwIn{$M^{-1}_{bjac},M^{-1}_{bjac}(p_l),r,p_w,p_l,relres,adp_{tol}$}
% 	\KwOut{$z$}
% 	\label{Algorithm 3.3}
	
% 	\eIf{ $relres \geq adp_{tol} $}{
		
% 		$z = \text{fMP-BJAC(}M^{-1}_{bjac}(p_l),r,p_w,p_l \text{)} $
% 	}{
% 		$ z = M^{-1}_{bjac}r$ 
% 	}
% 	return $z$
	
% \end{algorithm}

\begin{algorithm}
\caption{aMP-BJAC(lh) preconditioner}
\label{Algorithm 3.3}
\begin{algorithmic}[1] % [1] 表示行号
\Require $M^{-1}_{bjac}, M^{-1}_{bjac}(p_l), r, p_w, p_l, relres, adp_{tol}$
% \Ensure $z$
\If{$relres \geq adp_{tol}$}
    \State $z = \text{fMP-BJAC}(M^{-1}_{bjac}(p_l), r, p_w, p_l)$
\Else
    \State $z = M^{-1}_{bjac} r$
\EndIf
\State \Return $z$
\end{algorithmic}
\end{algorithm}

\section{Performance Evualation and Feature Analysis}\label{sec:4}

The proposed mixed precision preconditioners have been integrated into the JXPAMG solver \cite{JXPAMG}. In this section, using two typical test problems, we evaluate the performance improvement and investigate the convergence behavior of the mixed precision BJAC preconditioners. 

\subsection{Test Problems}\label{sec:4.1}

We consider two test problems: the three-dimensional diffusion equation (Diff3D) which is a model problem and the three-dimensional radiation hydrodynamics equation (RHD3D) which is a practical problem. For the Diff3D case, four different diffusion coefficients are concerned. For the RHD3D cases, two typical physical modelings are used. The detailed information of these test problems are shown in Table \ref{tab4.1}.

\begin{table*}[htbp]
	
	\centering
	\caption{Test problems.}
	\label{tab4.1}
	
	% \begin{tabular*}{0.91\linewidth}{lllccc} \toprule
 \begin{tabular*}{\linewidth}{llcccc} \toprule
		 Problem  &  Coefficients or Modeling  & Discretization &  Mesh Size & $N$ & $nnz$ \\ 
		\hline
		
\rule{0pt}{10pt} 		Diff3D-Const & isotropic & \multirow{4}{*}{7-point FDM } & \multirow{4}{*}{$128^3$} & \multirow{4}{*}{2,097,152} & \multirow{4}{*}{14,581,760} \\ 
		
\rule{0pt}{10pt} 		Diff3D-Ani($s$)   & anisotropic   &   &  & & \\ 
		
\rule{0pt}{10pt} 		Diff3D-Dis($s$)   & discontinuous &   &  & & \\ 
		
\rule{0pt}{10pt} 		Diff3D-Rand($s$)  & randomised    &   &  & & \\ 
		\hline
		
\rule{0pt}{10pt} 		RHD3D-1T  & single-temperature & \multirow{2}{*}{7-point FVM } & \multirow{2}{*}{$128^3$}  & 2,097,152 & 14,581,760   \\   
		
\rule{0pt}{10pt} 		RHD3D-3T   & three-temperature &   &  & 6,291,456  & 52,133,888 \\       
		\bottomrule
	\end{tabular*}

     \footnotesize $s$ in Diff3D-Ani(s), Diff3D-Dis(s), Diff3D-Rand(s) is the strength of the anisotropic, discontinuous and random coefficients, as described in Section \ref{sec:4.1.1}; FDM and FVM denote Finite Difference and Finite Volume Method respectively; $N$ and $nnz$ respectively denote the size and the number of non-zero entries of the matrix.
     
\end{table*}

\subsubsection{3D Diffusion Equation (Diff3D)}\label{sec:4.1.1}

The 3D diffusion equation with the Dirichlet boundary condition is given as follows:

\begin{equation}
	\label{eq4-1}
	\begin{cases}
		\begin{aligned}
			-\nabla (\kappa \nabla u ) &= f \text{,\quad } x \in \Omega \\
			u &= 0 \text{,\quad } x \in \partial \Omega \\
		\end{aligned}
	\end{cases}
\end{equation}
where $\Omega = (0,1)^3$ , $\kappa$ is the diffusion coefficient, which we consider the following four cases: 

\begin{enumerate}
	\item Constant case (Diff3D-Const): $\kappa = 1 $ in $\Omega$.
	
	\item Anisotropic case (Diff3D-Ani($s$)):

		% \begin{math}
            $$
		\kappa =
		% \left[
		\begin{pmatrix}
			1 & 0 & 0 \\
			0 & s & 0 \\
			0 & 0 & s \\
		\end{pmatrix}
  \text{in } \Omega \text{,}
		% \right]
  $$
	% \end{math}

where $s>1$ is the strength of anisotropic.
	
	\item Discontinuous case (Diff3D-Dis($s$)): 
	\begin{align*}
		\begin{split}
			\kappa = \left \{
			\begin{array}{ll}
				s,  & x \in {[0.25,0.75]}^3\\
				1,  & otherwise\\
			\end{array}, 
			\right.
		\end{split}
	\end{align*}
	where $s>1$ is the strength of discontinuous.
	
	\item Random case (Diff3D-Rand($s$)):  $\kappa = s^{\delta}$, where $0 \leq \delta \leq 1$ is a random function and $s>1$ is the strength of random.
	
\end{enumerate}

As shown in Table \ref{tab4.1}, these four test problems use the same discretization method. The mesh size is  ${128}^3$, and the size of the resulting matrix is 2,097,152.

\subsubsection{3D Radiation Hydrodynamic Equations (RHD3D)}\label{sec:4.1.2}

The radiation hydrodynamics equation (RHD) is a fundamental governing equation in the field of high energy density physics such as laser fusion and astrophysics \cite{Baldwin, Mo_2D3T, Xu_2D3T}. The operator splitting method is commonly used in practical simulations to divide this equation into hydrodynamic and radiation-diffusion equations to be solved separately, where the radiation-diffusion equation is shown in (\ref{eq4-2}).

\begin{equation}
	\label{eq4-2}
 \resizebox{0.95\hsize}{!}{$
	\begin{cases}
		\begin{aligned}
			c_{vr} \frac{\partial T_r}{\partial t} - \frac{1}{\rho} \nabla (K_r \nabla T_r ) &= \omega_{er} (T_e -T_r) \\
			c_{ve} \frac{\partial T_e}{\partial t} - \frac{1}{\rho} \nabla (K_e \nabla T_e ) &= \omega_{ei} (T_i -T_e) + \omega_{er} (T_r -T_e) \\
			c_{vi} \frac{\partial T_i}{\partial t} - \frac{1}{\rho} \nabla (K_i \nabla T_i ) &= \omega_{ei} (T_e -T_i) \\
		\end{aligned}
	\end{cases}
 $}
\end{equation}
where $\rho$ represents the medium density, $T_r$, $T_e$, and $T_i$ indicate the temperatures of photons, electrons, and ions, respectively; $c_{vr}$, $c_{ve}$, and $c_{vi}$ denote the isovolumetric specific heats of photons, electrons, and ions, respectively; $K_r$, $K_e$, and $K_i$ denote the diffusion coefficients; $\omega_{ei}$ and $ \omega_{er}$ denote the energy exchange coefficients of electron-ion and electron-photon, respectively.

Two modelings are concerned as follows:

\begin{enumerate}
	
	\item Three-temperatures modeling (RHD3D-3T). The resulting matrix for discreted RHD equation in (\ref{eq4-2}) has the following form in (\ref{eq4-3}):
	
	\begin{equation}
		\label{eq4-3}
		A=\left(
		\begin{matrix}
			A_R & D_{RE} & 0 \\
			D_{ER} & A_{E} & D_{EI} \\
			0 & D_{IE} & A_I
		\end{matrix} 
		\right) 
	\end{equation}
	where $A_R$, $A_E$, and $A_I$ are discrete systems of the three temperatures equations in (\ref{eq4-2}) with the same sparse pattern, and the matrices $D_{RE}$, $D_{RE}$, and $D_{EI}$, $D_{IE}$ are diagonal matrices reflecting the coupling of the three temperatures in (\ref{eq4-2}). 
	
	\item Single-temperature modeling(RHD3D-1T): if the three-temperatures reach to equilibrium state, the equation degenerates into a scalar single-temperature equation with the matrix of $A_R$ in (\ref{eq4-3}).
	
\end{enumerate}

It should be pointed out that the RHD3D1T and RHD3D3T are typical multiscale systems that are challenging to solve and have been chosen as the SolverChallenge competition problems (\url{https://www.solver-conference.cn/solverchallenge23/index.html}).

\subsection{Test Setting}\label{sec:4.2}

For all experiments in this paper, the BJAC and its mixed precision versions fMP-BJAC, aMP-BJAC(hl) and aMP-BJAC(lh) are all used as the preconditioners for PCG. A reduction tolerance of relative residual L2-norm(RelResNorm) $tol_{res}$ is used as the convergence criterion for PCG iteration. $tol_{res}$ is chosen depend on the problem, for the model problem Diff3D, $tol_{res}=10^{-10}$, while for the practical problem RHD3D, $tol_{res}=10^{-12}$. The zero vector is used as the initial guess solution of iterations. 

For two kinds of test problems, single precision (fp32) is used as the low precision format. For working precidion (high precision) format, double precision (fp64) is used for the model problem Diff3D, and long double precision (fp80) is used for the practical problem RHD3D. 

All experiments are performed on a machine with multicore compute nodes, the number of processes and the number of blocks in BJAC both equal to 32. Unless otherwise specified, the iteration numbers of inter-block and intra-block in the BJAC are set to $k=t=2$.

We denote the following notations used in the numerical results in the following section.

\begin{itemize}
\item ``fp32'', ``fp64'', ``fp80'' denote single, double, and long double precision respectively.
\item ``fp64-uniform-BJAC PCG'' and ``fp80-uniform-BJAC PCG'' respectively denote the uniform precision using fp64 and fp80 as working precision for BJAC preconditioner and PCG iteration.
\item  ``fp32-fMP-BJAC PCG'', ``fp32-aMP-BJAC(hl) PCG'' and ``fp32-aMP-BJAC(lh) PCG'' respectively denote the fMP-BJAC, aMP-BJAC(hl), and aMP-BJAC(lh) with fp32 for BJAC preconditioners and working precision(high precision) for PCG.
\end{itemize}

\subsection{Performance Evaluation}\label{sec:4.3}

In this section, we evaluate the performance gain of the proposed mixed precision preconditioners. Define the speedup as follows:

\begin{equation}
\label{eq4-4}
    Speedup = \frac{T_{uniform}}{T_{mix}}
\end{equation}
where $T_{uniform}$ denotes the CPU time of the uniform precision BJAC PCG algorithm and $T_{mix}$ denotes the CPU time of the mixed precision BJAC PCG algorithm, i.e., fMP-BJAC or aMP-BJAC PCG algorithm. 
For aMP-BJAC(hl) and aMP-BJAC(lh), the performance is depend on the adaptive threshold ${adp}_{tol}$.  
In the experiments in this section, unless otherwise statement, $adp_{tol}=10$ for aMP-BJAC(hl), and  $adp_{tol}={10}^{-5}$ for aMP-BJAC(lh).

\begin{figure*}[ht]
    \centering
    \includegraphics[width=0.88\linewidth]{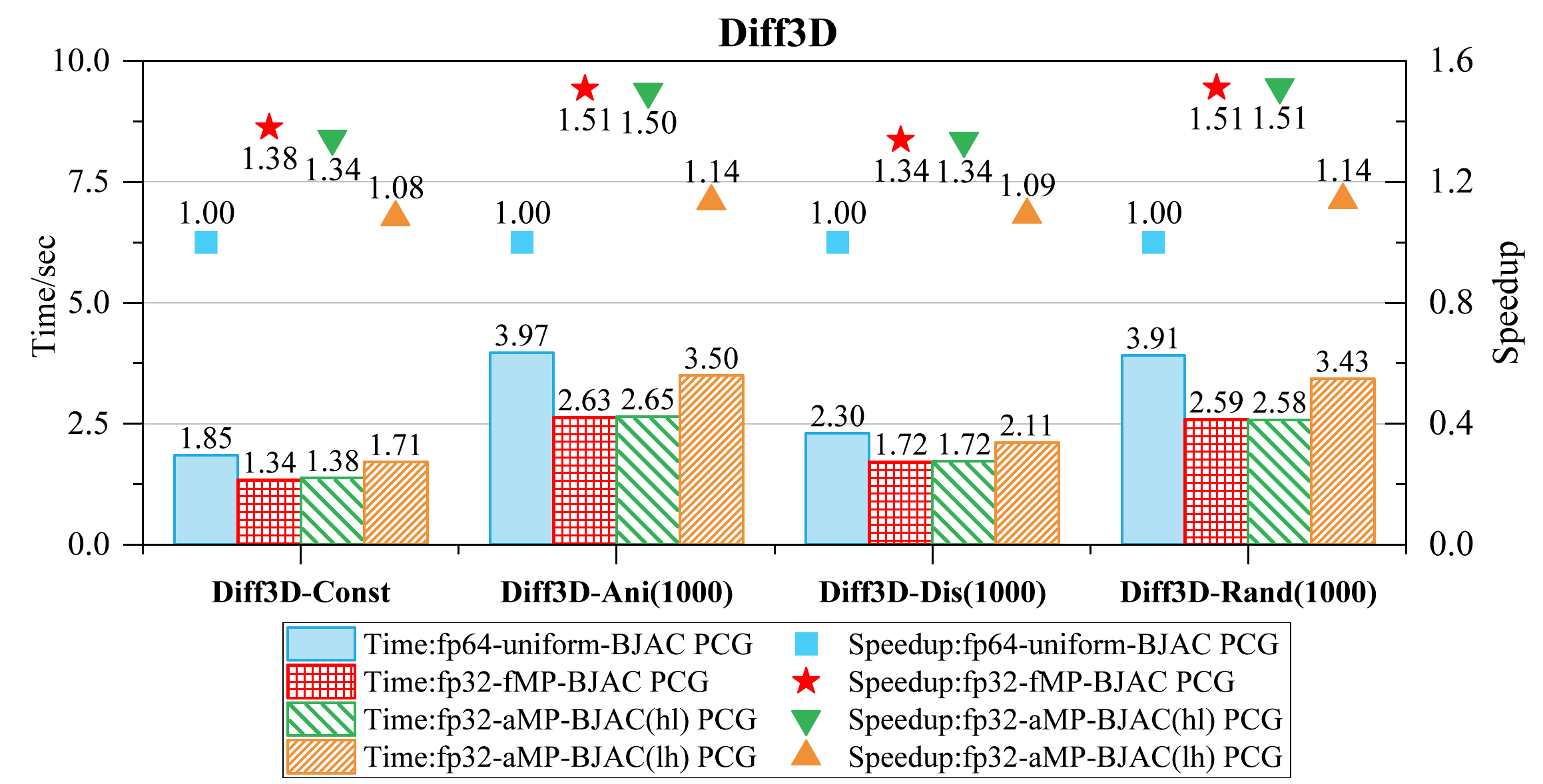}
    \caption{Time and speedup of mixed precision BJAC PCG compared to fp64-uniform-BJAC PCG: Diff3D problems.}
    \label{fig4.1}
\end{figure*}

\begin{figure*}[ht]
    \centering
    \includegraphics[width=0.88\linewidth]{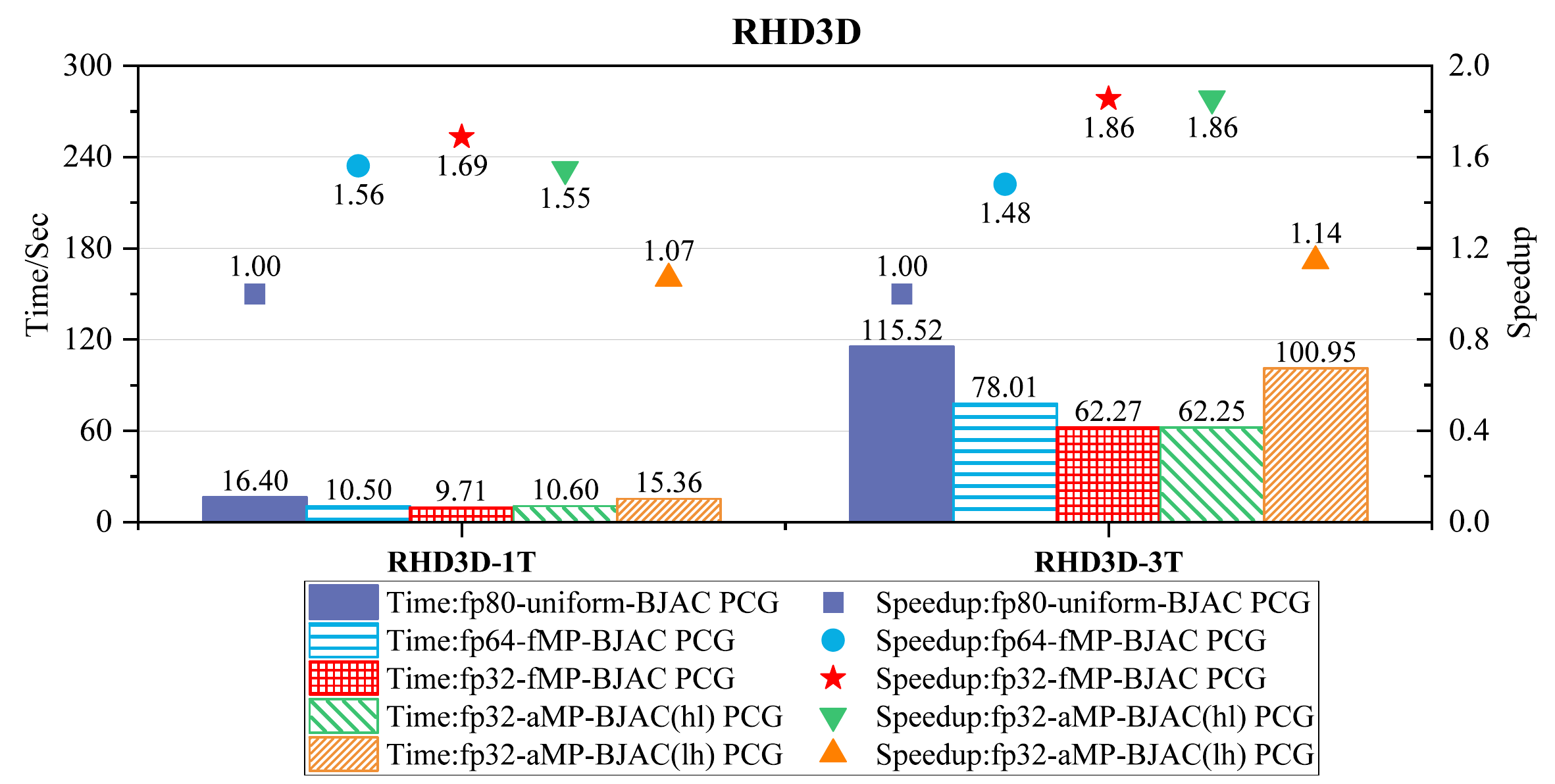}
    \caption{Time and speedup of mixed precision BJAC PCG compared to fp80-uniform-BJAC PCG: RHD3D problems.}
    \label{fig4.2}
\end{figure*}

The results for Diff3D and RHD3D problems are given in 
Figure \ref{fig4.1} and Figure \ref{fig4.2} respectively. We can conclude from the results that, for all test problems, the proposed three mixed precision preconditioners can improve performance in CPU time compared to the uniform precision BJAC preconditioner.

For the Diff3D problems, fp64 is chosen as the high precision and fp32 is chosen as the low precision. From Figure \ref{fig4.1}, we can see that different mixed precision preconditioner have different performance behavior, and the performance gains differ from diffusion coefficients. For all cases, fMP-BJAC and aMP-BJAC(hl) have achieved comparable speedup, that of 1.3x to 1.5x, which outperform the speedup of aMP-BJAC(lh). 

For the RHD3D problems, fp80 is chosen for the high precision and fp32 for the low precision. From Figure \ref{fig4.2}, we can conclude that the fMP-BJAC and aMP-BJAC(hl) can achieve speedup of more than 1.5$\times$ for RHD3D-1T and 1.8$\times$ for RHD3D-3T respectively, and fMP-BJAC slightly outperforms aMP-BJAC(hl) for the RHD3D-1T case. Again, the speedups of both fMP-BJAC and aMP-BJAC(hl) are better than that of aMP-BJAC(lh) for two cases. 
In addition, we also investigate the results of fMP-BJAC with fp64 as the low precision (fp64-fMP-BJAC PCG in Figure \ref{fig4.2}), its speedup does not outperform that of fp32-fMP-BJAC PCG for both cases.

From the above results, it can be concluded that the mixed precision preconditioners fMP-BJAC and aMP-BJAC(hl) achieve performance gains compared to the uniform precision BJAC. Meanwhile, two issues should be noticed. Firstly, the adaptive precision with two orders, i.e., aMP-BJAC(hl) and aMP-BJAC(lh), have quite performance differences, and the "high-to-low" strategy is significantly better than the "low-to-high" strategy. Secondly, the mixed precision preconditioners have different performance gains for different problems. The first issue concerns the influence of the mixed precision strategy on the convergence behavior, and the second issue involves the impact of the application features on the performance of the mixed precision preconditioner. We furtherly analyse and discuss these two issues in the following section.

\subsection{Convergence Behavior Analysis}\label{sec:4.4}

One of the main factors that impact the performance of mixed preconditioners is their  convergence behavior. Actually, the key to achieving maximum performance speedup is the trade-off between the computational cost of single iteration and the
convergence speed of iterations. In this section, we analyse the convergence behavior of the proposed mixed precision BJAC preconditioners.

\subsubsection{Convergence Delays Phenomenon}\label{sec:4.3.1}

Ideally, it is expected that the convergence speed of the preconditioner, which measured by iteration numbers that reach to the convergence criterion, does not deteriorate due to the mixed precision. However, in some problems, the numerical results show that, compared to the uniform precision BJAC, additional iterations are required for the mixed precision BJAC, which we called the convergence delay phenomenon. 

\begin{figure*}[ht!]
    \centering
    \includegraphics[width=0.7\linewidth]{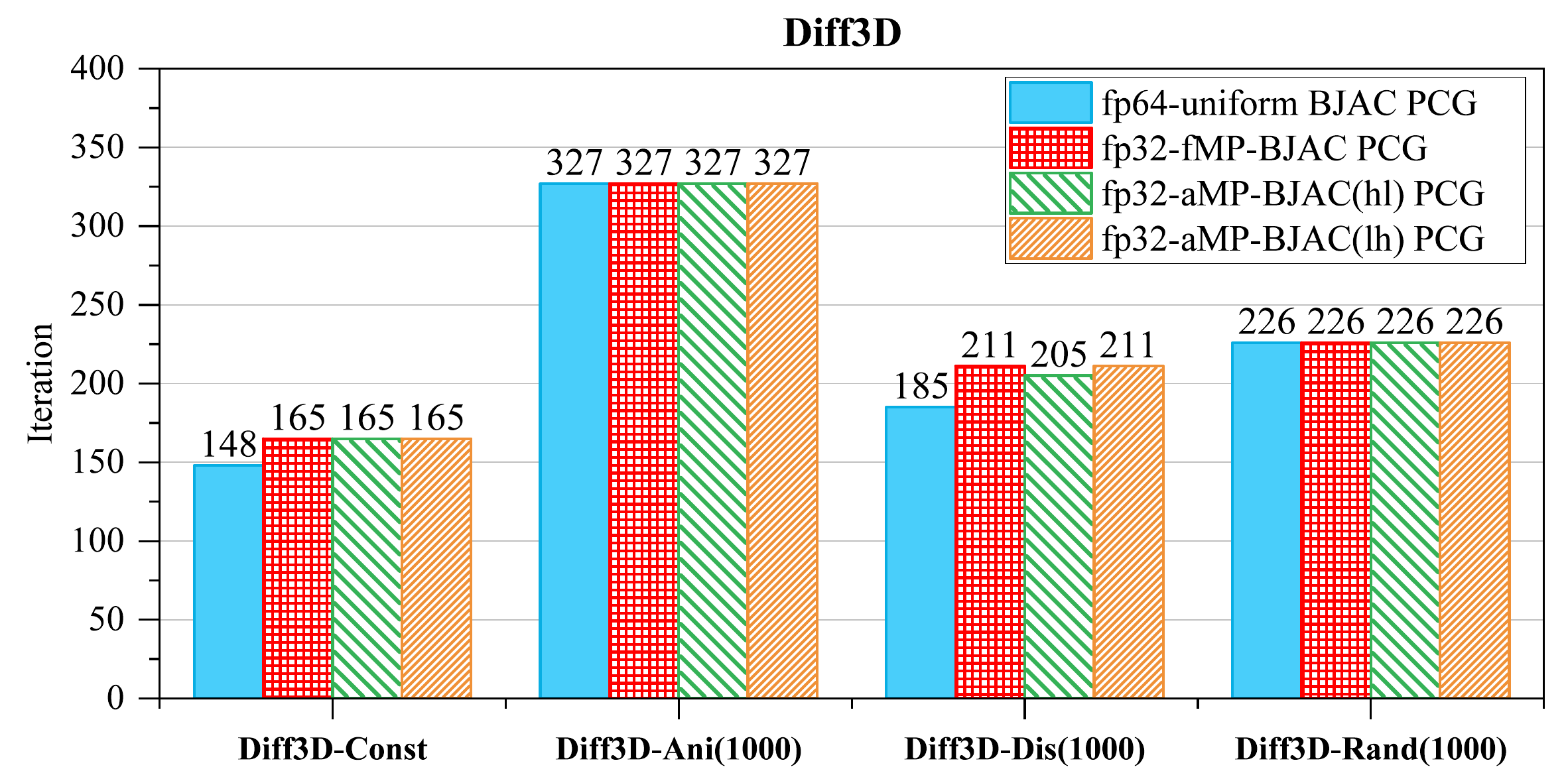}
    \caption{Iteration numbers for fp64-uniform- and mixed precision BJAC: Diff3D problems.}
    \label{fig4.3}
\end{figure*}

\begin{figure*}[ht!]
    \centering

    \includegraphics[width=0.7\linewidth]{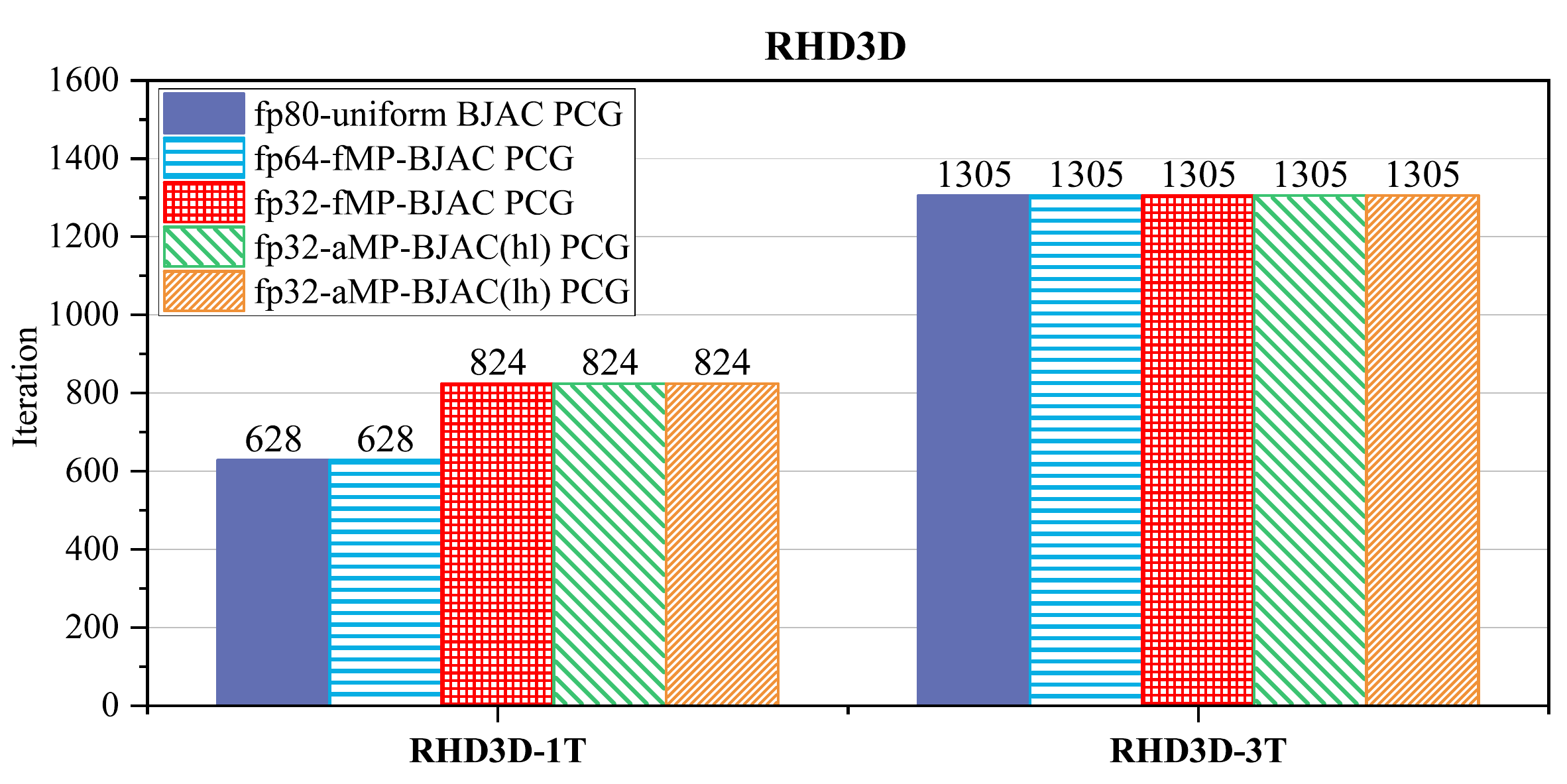}
        \caption{Iteration numbers for fp80-uniform- and mixed precision BJAC: RHD3D problems.}
        \label{fig4.4}
\end{figure*}

Corresponding to the tests in Figures \ref{fig4.1} and \ref{fig4.2}, the results of iteration numbers for the Diff3D and RHD3D problems are given in Figures \ref{fig4.3} and \ref{fig4.4} respectively.

From the results of Diff3D problems in Figure \ref{fig4.3}, it can be observed that the convergence behavior of mixed precision preconditioners varies with diffusion coefficients. For Diff3D-Ani(1000) and Diff3D-Rand(1000) cases, the mixed precision BJAC have the same iteration numbers as the fp64-uniform-BJAC, which achieve the perfect performance. However, for other two problems, Diff3D-Const and Diff3D-Dis(1000), the iteration numbers of the mixed precision preconditioners increases slightly compared to the fp64-uniform-BJAC. Concretely, for Diff3D-Const, the three mixed precision preconditioners increase iteration numbers by 11\%, which is from 148 to 165. For Diff3D-Dis(1000), the convergence behavior is different due to preconditioners, fp32 fMP-BJAC and aMP32 MP-BJAC(lh)  increase the iteration numbers by 14\% which is from 185 to 211, and the fp32 aMP-BJAC(hl) increases the iteration numbers by 11\% which is from 185 to 205.

For RHD3D problems, the results in Figure \ref{fig4.4} show that the convergence behavior of mixed precision preconditioners is different due to physical modeling. For RHD3D-3T case, all mixed precision preconditioners achieve perfect convergence, i.e., iteration numbers are the same as that of the f80-uniform-BJAC. However, for RHD3D-1T case, compared to the fp80-uniform-BJAC, the three mixed precision BJAC preconditioners which used fp32 as low precision, include fp32-fMP-BJAC, fp32-aMP-BJAC(hl)
and fp32-aMP-BJAC(lh), increase the number of iterations by 31\% which from 628 to 824. 
It can be noticed from Figure \ref{fig4.4} that, for the RHD3D-1T problem, if fp64 is used as low precision, the mixed precision preconditioner fp64-fMP-BJAC can achieve the same iteration numbers as fp80-uniform-BJAC preconidioner.

In summary, for three test cases concerned in this paper, include Diff3D-Const, Diff3D-Dis(1000) and RHD3D-1T, the mixed precision BJAC preconditioners lead to the convergence delays issue, additional iterations are required to achieve the convergence criterion compared to the uniform-BJAC preconditioner. This issue is crucial for achieving better speedup of mixed precision preconditioners, and will be discussed more detail in the following section.

%% ---------------------------------------------------------------------------

\subsubsection{Improving Convergence via Adaptive Precision Preconditioners}\label{sec:4.3.2}

For adaptive mixed precision preconditioners aMP-BJAC(hl) and aMP-BJAC(lh), the convergence behavior is influenced by the adaptive threshold ${adp}_{tol}$. Hence we further investigate the impact of adaptive threshold on the convergence behavior of these two preconditioners.

In general, given a threshold ${adp}_{tol}$, the aMP-BJAC(hl) preconditioner firstly executes high precision iterations and then executes low precision in the rest of iterations, while the aMP-BJAC(lh) executes iterations in the opposite order. So the proportion and the distribution of the low precision executed during the iterations can be manipulated by choosing ${adp}_{tol}$. For aMP-BJAC(hl), generally, the larger ${adp}_{tol}$ is,  the larger the proportion of low precision iterations is. In extreme case, the aMP-BJAC(hl) degenerates into the fMP-BJAC when ${adp}_{tol}$ is sufficiently large, otherwise, it degenerates into a uniform-BJAC preconditioner when ${adp}_{tol}$ is sufficiently small. For the aMP-BJAC(lh), 
contrarily in general, the larger ${adp}_{tol}$ is, the smaller the proportion of low precision iterations is. In extreme case, the aMP-BJAC(lh) respectively degenerates into the uniform-BJAC and fMP-BJAC preconditioner when ${adp}_{tol}$ is sufficiently large and small.

For three test cases, include Diff3D-Const, Diff3D-Dis(1000) and RHD3D-1T, where the mixed precision BJAC preconditioners cause convergence delays as shown in last section, Table \ref{tab4.2} gives the results of convergence behavior, i.e., changing trend of iteration numbers, by varying the threshold ${adp}_{tol}$ with typical values.

\begin{table}[ht]
    \centering
    \caption{The iteration numbers of aMP-BJAC(hl) and aMP-BJAC(lh) preconditioner varying with ${adp}_{tol}$ for three typical test problems.}
    \label{tab4.2}
    \begin{tabularx}{0.477\textwidth}{p{2.16cm}|p{0.64cm}|p{0.76cm}p{0.97cm}p{1.2cm}}
        \toprule
        \multirow{2}{*}{Preconditioner} &\multirow{2}{*}{${adp}_{tol}$}&  \multicolumn{3}{l}{Iteration Numbers of PCG}\\
        
          &    &    Diff3D-Const & Diff3D-Dis(1000) & RHD3D-1T\\
        \hline
        
        uniform-BJAC  & - & 148 & 185 & 628\\
        \hline
        
        fp32-fMP-BJAC & - & 165 & 211 & 824\\
        \hline
        
        \multirow{8}{*}{aMP-BJAC(hl) }   & $10$    & 165  & 205  & 824\\
                        & $5$     & 148  & 201  & 824\\
                        & $1$     & 148  & 185  & 824\\
                        & $10^{-1}$ & 148 & 185 & 628\\
                        & $10^{-2}$ & 148 & 185 & 628\\
                        & $10^{-3}$ & 148 & 185 & 628\\
                        & $10^{-4}$ & 148 & 185 & 628\\  
                        & $10^{-5}$ & 148 & 185 & 628\\ 
        \hline
        \multirow{8}{*}{aMP-BJAC(lh) }    & $10^{-5}$  & 165 &  211  & 824\\
                        & $10^{-4}$  & 165 &  211  & 824\\
                        & $10^{-3}$  & 165 &  211  & 824\\
                        & $10^{-2}$  & 165 &  211  & 824\\
                        & $10^{-1}$  & 165 &  211  & 824\\
                        & $1$      & 165 &  211  & 824\\
                        & $5$      & 165 &  211  & 824\\
                        & $10$     & 165 &  201  & 824\\
        \bottomrule      
    \end{tabularx}    
\end{table}

The results in Table \ref{tab4.2} show that, for these three test cases, as the adaptive threshold ${adp}_{tol}$ decreases, the iteration numbers of the aMP-BJAC(hl) decreases and eventually converge to that of uniform-BJAC. However, for aMP-BJAC(lh), as the ${adp}_{tol}$ increases, the iteration numbers do not converge to that of uniform-BJAC, even do not decrease for Diff3D-Const and RHD3D-1T cases. It indicates the notable differences in convergence behavior between these two preconditioners.

\begin{figure}[ht]
    \centering

    \includegraphics[width=1\linewidth]{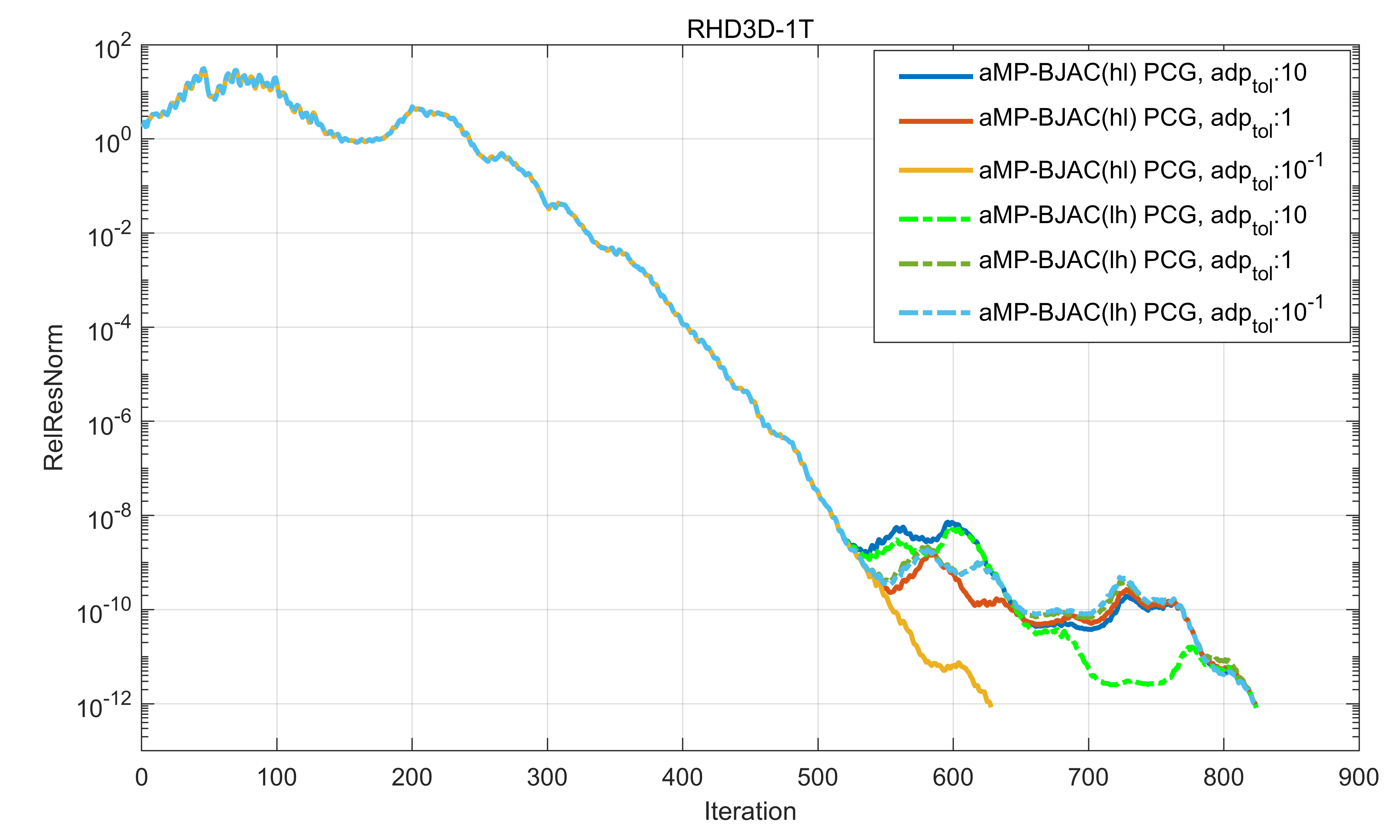}\\
    \includegraphics[width=0.96\linewidth]{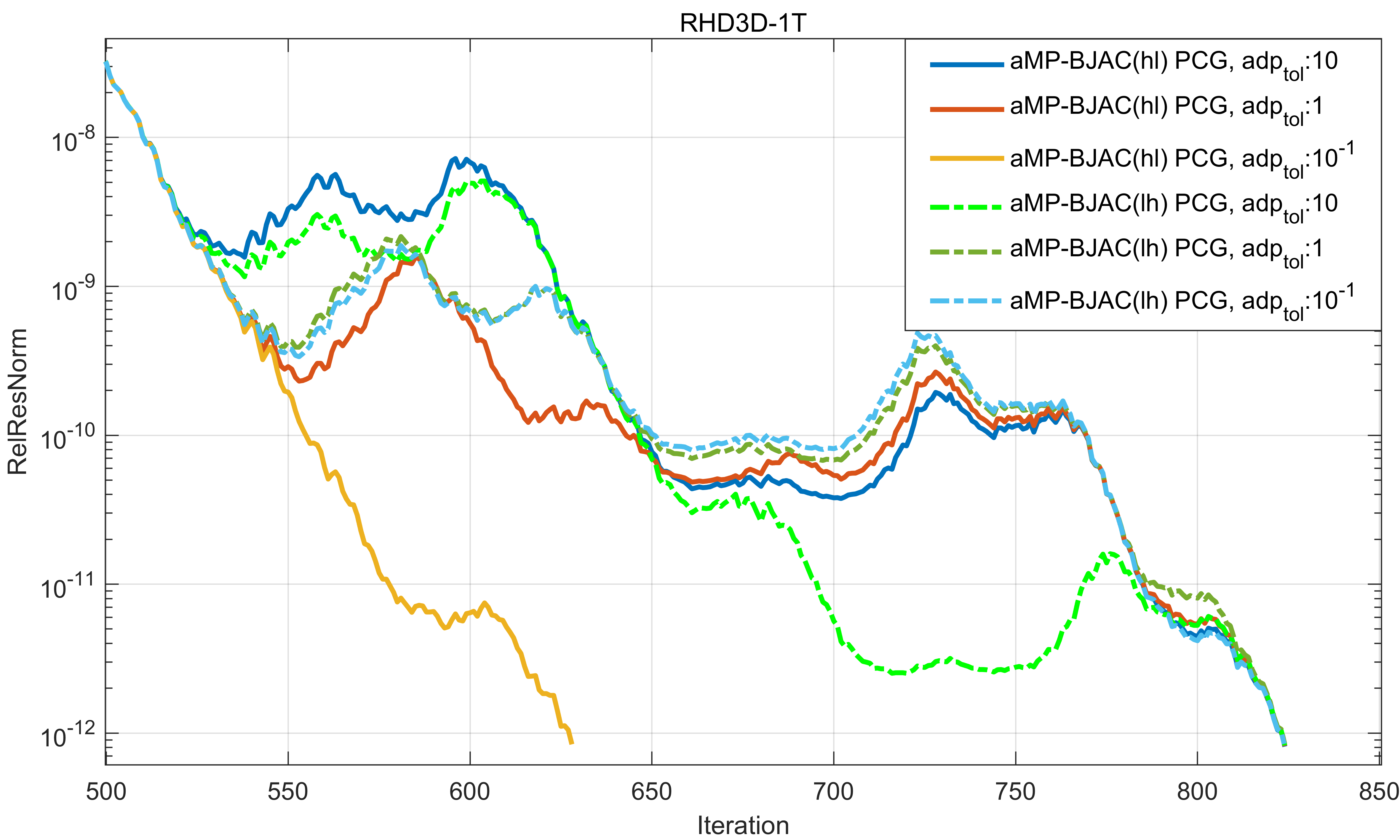}
        
    \caption{RelResNorm curves of the aMP-BJAC(hl) and aMP-BJAC(lh) for solving the RHD3D-1T problem with three typical values of ${adp}_{tol}$. Upper panel: the overall result. Lower panel: the detailed result after the 500th iteration.}
    \label{fig4.5}
\end{figure}

Furthermore, taking the RHD3D-1T as an example, we further investigate the reduction behaviour of the relative residual L2-Norm (RelResNorm) for the aMP-BJAC(lh) and aMP-BJAC(hl). The results of these adaptive mixed precision preconditioners with three typical values of adaptive thresholds ${adp}_{tol} = 10, 1, {10}^{-1}$ are given in Figure \ref{fig4.5} 

As shown in Figure \ref{fig4.5}, in the early and middle stages during the iterations, the RelResNorm curves of six preconditioners exactly overlap. When the RelResNorm is reduced to around $5\times {10}^{-10}$, only the RelResNorm of the aMP-BJAC(hl) with ${adp}_{tol} = {10}^{-1}$ continue to reduce, however, for other values of ${adp}_{tol}$, the RelResNorm curves occur bifurcation and reach the convergence criterion with additional iterations, i.e., the convergence delays occur. In particular, for the aMP-BJAC(lh), even if ${adp}_{tol}$ is increased to $10$, at which only the beginning 1st, the 57th-60th, and the 62nd-102nd iterations use low precision, and the rest of the iterations use high precision, the convergence delays still occur.

The analyses above show that the aMP-BJAC(hl) preconditoner with a suitable adaptive threshold can improve or even avoid the the convergence delay of the fMP-BJAC preconditioner. 
Meanwhile, as shown in Table \ref{tab4.2}, the convergence behavior of aMP-BJAC(lh) preconditioner indicates that it is hard to avoid the convergence delay if the iteration starts with the low precision for even just few iterations, which needs further theoretical analyses in the future.

\subsection{Feature Analysis for Convergence Delays}\label{sec:4.5}

Convergence delay leads to a challenge for achieving ideal speedup of mixed preconditioners. 
An interesting observation is that, the three problems occurring the convergence delay, include Diff-Const, Diff-Dis(1000), and RHD3D-1T, are more closer to a model problem compared to other three test problems, and should intuitively more likely to achieve the ideal speedup, whereas the numerical results are just the opposite. 
This is an interesting phenomenon that should be concerned about in the future.

In this section, we investigate the correlation between the convergence delays behaviour and the problem features. As shown in Table \ref{tab4.1}, the four Diff3D test cases and the RHD3D-1T problem have the same matrix sparse pattern as well as bandwidth since all tests use the 3D 7-point scheme as discretization. For the RHD3D-3T test problem, using the same discretization mesh and scheme as other test problems, however, its matrix size is three times larger than other test problems since the coupling of three-temperatures involved on each mesh cell. In addition, for all test problems, the resulting matrices are belong to M-matrix\cite{Saad}.

Therefore, the sparsity patterns and numerical properties of these matrices are overall similar, while the main differences among these matrices are magnitudes in entries value of matrix. Concretely, we consider two features that related to the magnitude differences of entries, the multiscale and diagonal dominance. We analyse the impact of these two features on the convergence delay caused by the mixed precision preconditioner.

\subsubsection{Multiscale Feature}\label{sec:4.5.1}

The multiscale of a matrix refers to the magnitude differences between the off-diagonal entries in the same row of the matrix. The definition and detailed discussion of multiscale matrices can refer to \cite{Xu_multiscale}, here we only give a basic concept of the multiscale matrix. For the matrix $A={(a_{ij})}_{n\times n} \in R^{n\times n}$, the multiscale strength measure $\tau_i$ of the $i$-th row is defined as follows:

\begin{equation}
    \label{eq4-5}
    \tau_i = \frac{\max \{ \vert a_{ij} \vert  \text{,} i \neq j \} }{\min \{ \vert a_{ij} \vert \text{,} i \neq j \} } , \text{ \; when \;} a_{ij} \neq 0.
\end{equation}
Given a strong multiscale threshold $\theta$, a matrix is called to be strongly multiscale if there exists a row $i$ whose multiscale strength is greater than $\theta$, i.e., $\tau_i \geq \theta$. Otherwise, if the multiscale strength of all rows is less than $\theta$, the matrix is weakly multiscale or single-scale. 

The multiscale property of matrices is a common feature in many practical applications. For the six problems considered in this paper, the multiscale features mainly arising from the non-smoothed diffusion coefficients, and coupling strength in RHD3D-3T case. It is clear that the constant coefficient problem Diff3D-Const is a single-scale matrix, while the other problems are multiscale matrices whose multiscale strength determined by $s$, which is the strength of discontinuity, anistropics, and random of the diffusion coefficients in Table \ref{tab4.1}. 

Table \ref{tab4.3} gives the distribution of the multiscale strength for the six problems.
The results illustrate that Diff3D-Const is a single-scale problem, while Diff3D-Dis(1000) and RHD3D-1T have relatively weak multiscale strength, and the other problems have strong multiscale. Associate with the results in Section \ref{sec:4.4} (see Figure \ref{fig4.3} and Figure \ref{fig4.4}), it can be found that the convergence delays of mixed precision preconditioners are related to the multiscale strength. Concretely, the cases occurring convergence delays are exactly single-scale or weak-multiscale problems. On the contrary, the strong multiscale problems, including Diff3D-Ani(1000), Diff3D-Rand(1000), and RHD3D3T, do not occur convergence delays.

\begin{table*}[!ht]
    \centering
    \caption{Multiscale strength distribution for Diff3D and RHD3D problems, i.e., the percentage of row numbers belong to different interval of multiscale strength.}
    \label{tab4.3}

    \begin{tabular}{c|rrrr|rr} \toprule
\rule{0pt}{10pt}         \multirow{2}{*}{\makecell[c]{Multiscale\\ Strength Interval}} & \multicolumn{4}{c|}{Diff3D} &  \multicolumn{2}{c}{RHD3D}\\
\rule{0pt}{10pt}           &  Const  &   Ani(1000) &   Dis(1000) &   Rand(1000)  &  1T &  3T\\
        \hline
\rule{0pt}{10pt}        $[1,\; \;  \; \quad   {10}^1)$ & 100\% & 0\% & 98.86\% & 49.01\% & 99.85\% & 8.24\% \\
\rule{0pt}{10pt}        $[{10}^1, \; \; \;    {10}^2)$ & 0\% & 0\% & 0\% & 40.90\% & 0.11\% & 4.58\% \\ 
\rule{0pt}{10pt}        $[{10}^2, \; \; \;    {10}^3)$ & 0\% & 0\% & 1.14\% & 10.09\% & 0.03\% & 8.94\% \\ 
\rule{0pt}{10pt}        $[{10}^3, \; \; \;    {10}^4)$ & 0\% & 100\% & 0\% & 0\% & 0.01\% & 28.87\% \\ 
\rule{0pt}{10pt}        $[{10}^4, \; \; \;    {10}^5)$ & 0\% & 0\% & 0\% & 0\% & 0\% & 8.93\% \\ 
\rule{0pt}{10pt}        $[{10}^5, \; \;    {10}^{10})$ & 0\% & 0\% & 0\% & 0\% & 0\% & 21.75\% \\ 
\rule{0pt}{10pt}        $[{10}^{10}, \:    {10}^{15})$ & 0\% & 0\% & 0\% & 0\% & 0\% & 7.35\% \\ 
\rule{0pt}{10pt}        $[{10}^{15}, \:   +\infty)$ & 0\% & 0\% & 0\% & 0\% & 0\% & 11.34\% \\
        \bottomrule
    \end{tabular}
\end{table*}

In the following, taking the multiscale problem Diff3D-Ani($s$) as an example, we investigate the correlation between the convergence delay and multiscale strength by tuning the parameter $s$. We consider six cases with $s=1000,100,10,4,2,1$, and the multiscale strength distributions of these six cases are given in Table \ref{tab4.4}. We can see that from Table \ref{tab4.4}, as $s$ decreasing, the multiscale strength become more weak. 

\begin{table*}[!ht]
    \centering
    \caption{Multiscale strength distribution for the Diff3D-Ani(s) with s=1000,100,10,4,2,1.}
    \label{tab4.4}
    \begin{tabular}{c|rrrrrr} \toprule
\rule{0pt}{10pt}         \multirow{2}{*}{\makecell[c]{  {Multiscale}\\ Strength Interval}} & \multicolumn{6}{c}{ Diff3D-Ani(s)} \\

\rule{0pt}{10pt}         & Ani(1000) &  Ani(100) & Ani(10) & Ani(4) & Ani(2) & Ani(1) \\ \hline

\rule{0pt}{10pt}        $[1, \qquad \; \; \;  \; \; 2)$ & 0\% & 0\% & 0\% & 0\% & 0\% & 100\% \\ 
\rule{0pt}{10pt}        $[2, \qquad \; \; \;  \; \; 4)$ & 0\% & 0\% & 0\% & 0\% & 100\% & 0\% \\ 
\rule{0pt}{10pt}        $[4,  \qquad \;  \;   \; 10)$ & 0\% & 0\% & 0\% & 100\% & 0\% & 0\% \\ 
\rule{0pt}{10pt}        $[10, \quad \; \; \;  100)$ & 0\% & 0\% & 100\% & 0\% & 0\% & 0\% \\ 
\rule{0pt}{10pt}        $[100,\; \;  \; 1000)$ & 0\% & 100\% & 0\% & 0\% & 0\% & 0\% \\ 
\rule{0pt}{10pt}        $[1000, \;  \; +\infty)$ & 100\% & 0\% & 0\% & 0\% & 0\% & 0\% \\  
        \bottomrule
    \end{tabular}
\end{table*}

\begin{figure*}[!ht]
    \centering
    \includegraphics[width=0.47\linewidth]{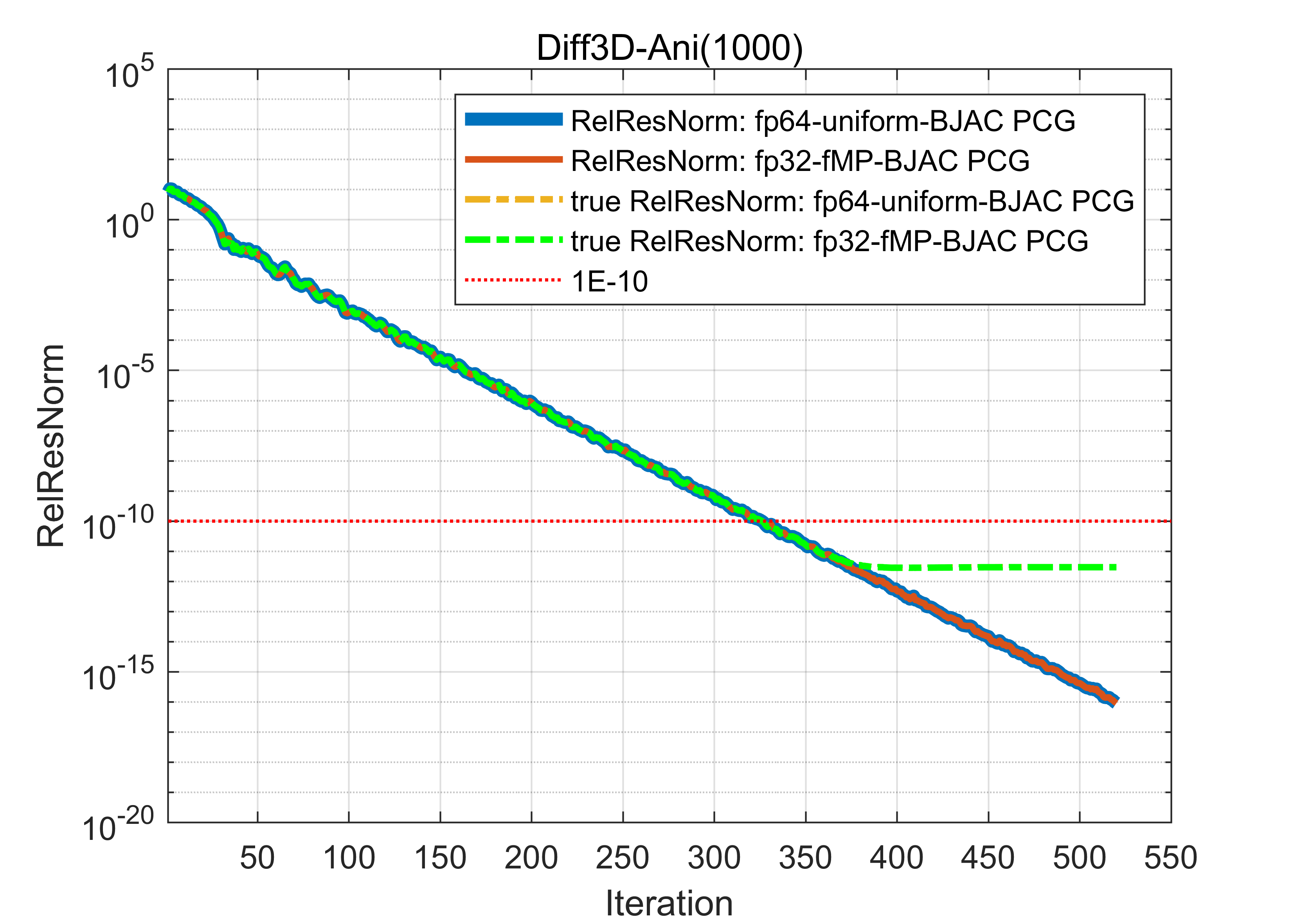}
    \includegraphics[width=0.47\linewidth]{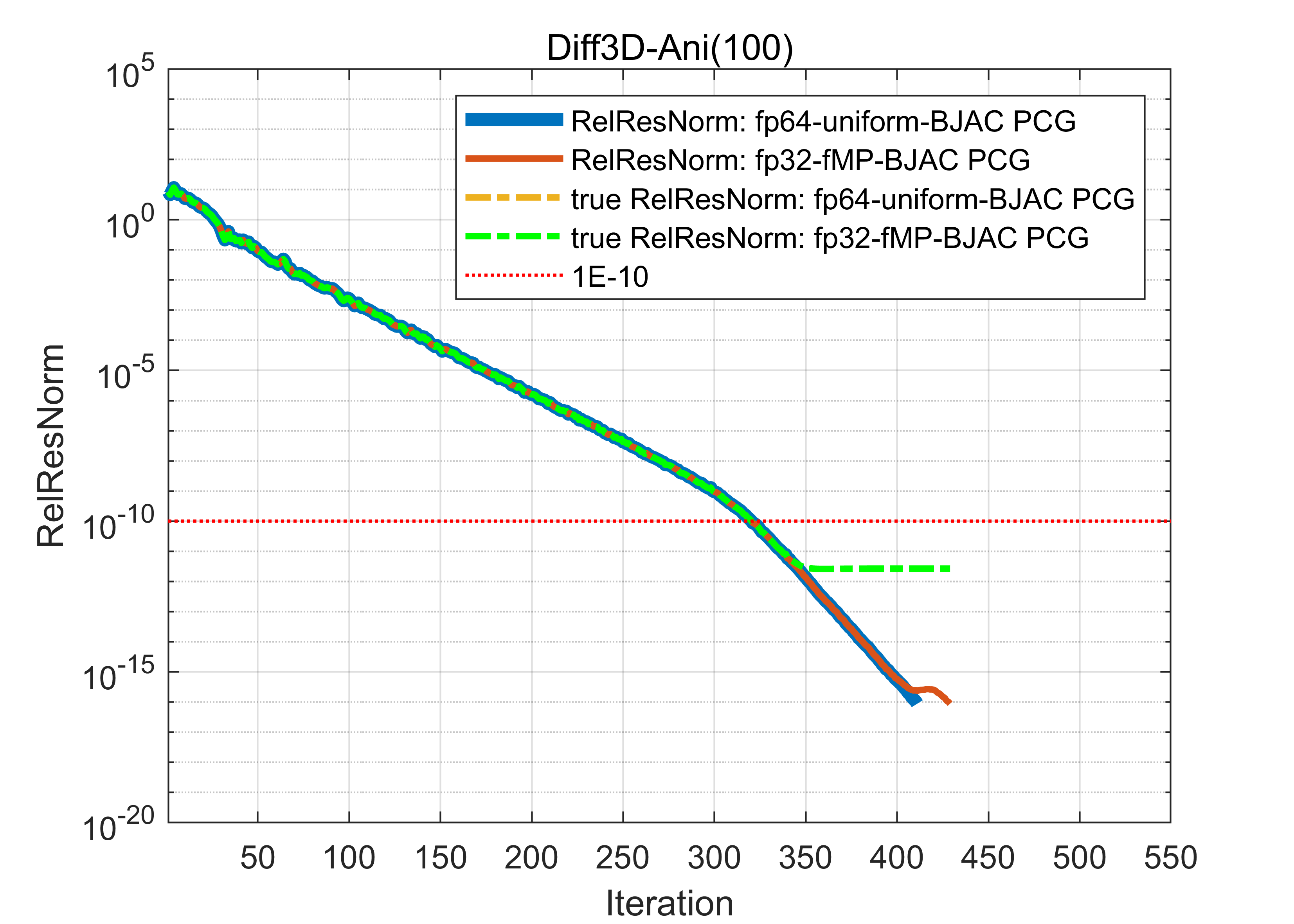} \\
    \includegraphics[width=0.47\linewidth]{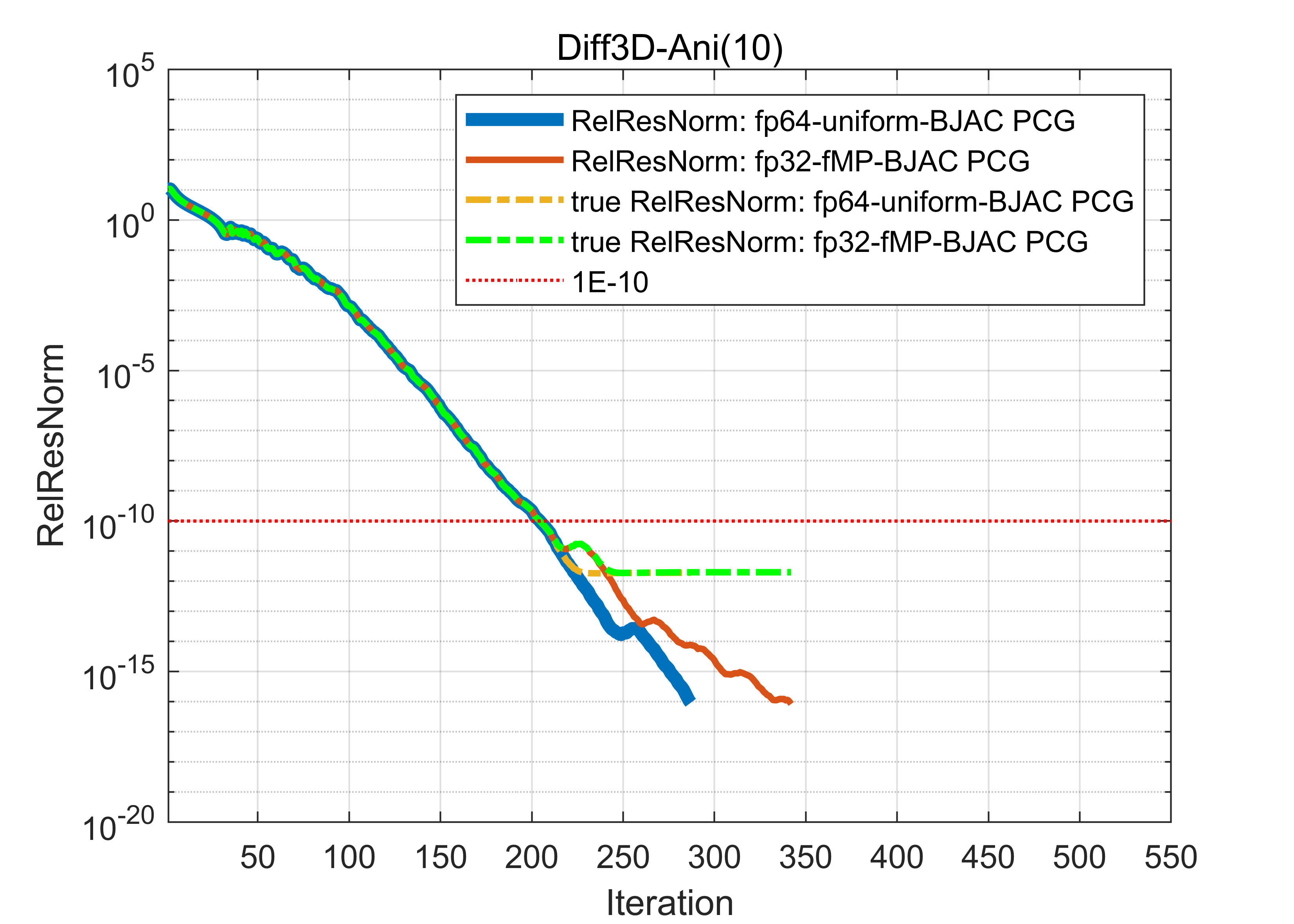}
    \includegraphics[width=0.47\linewidth]{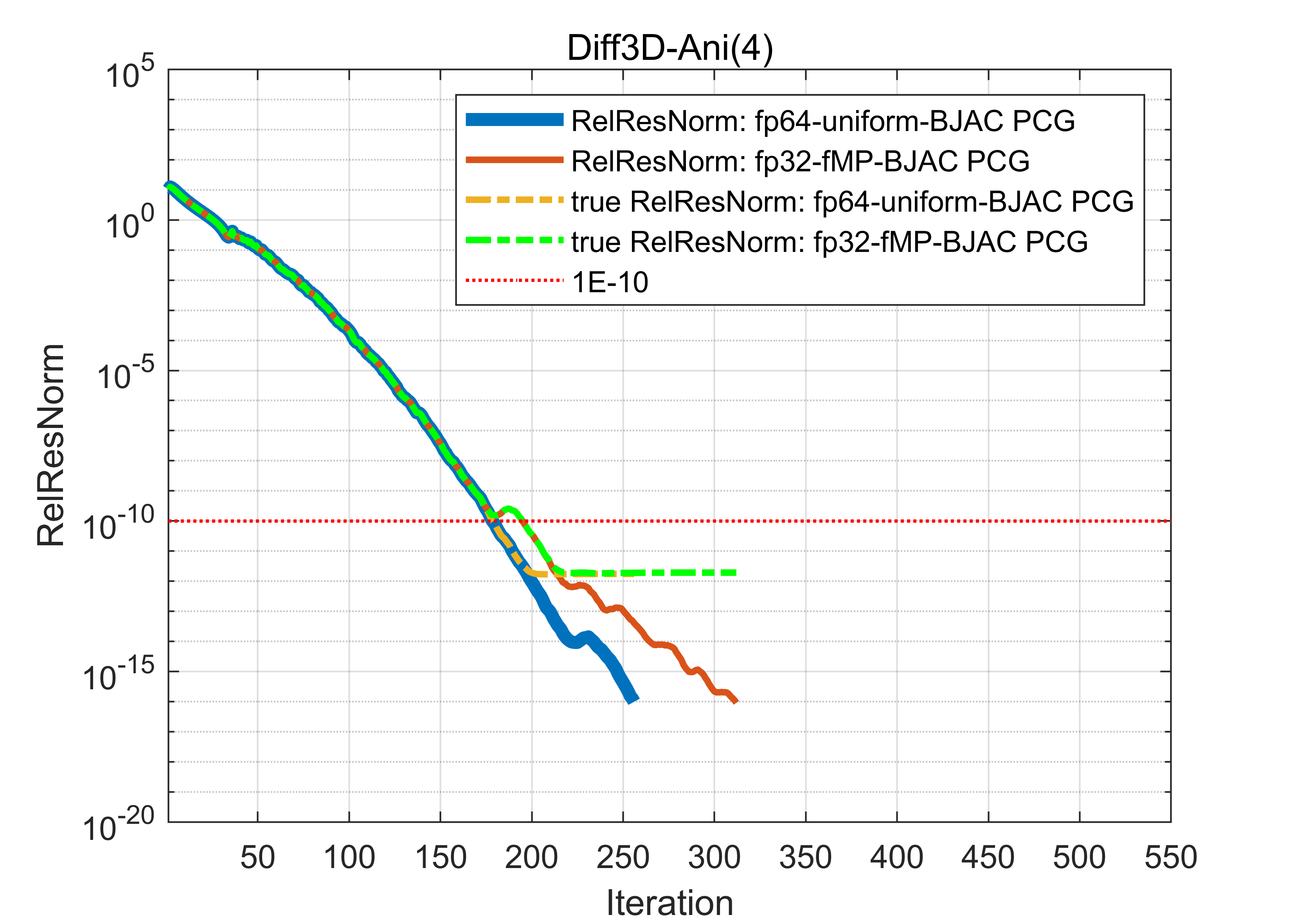} \\
    \includegraphics[width=0.47\linewidth]{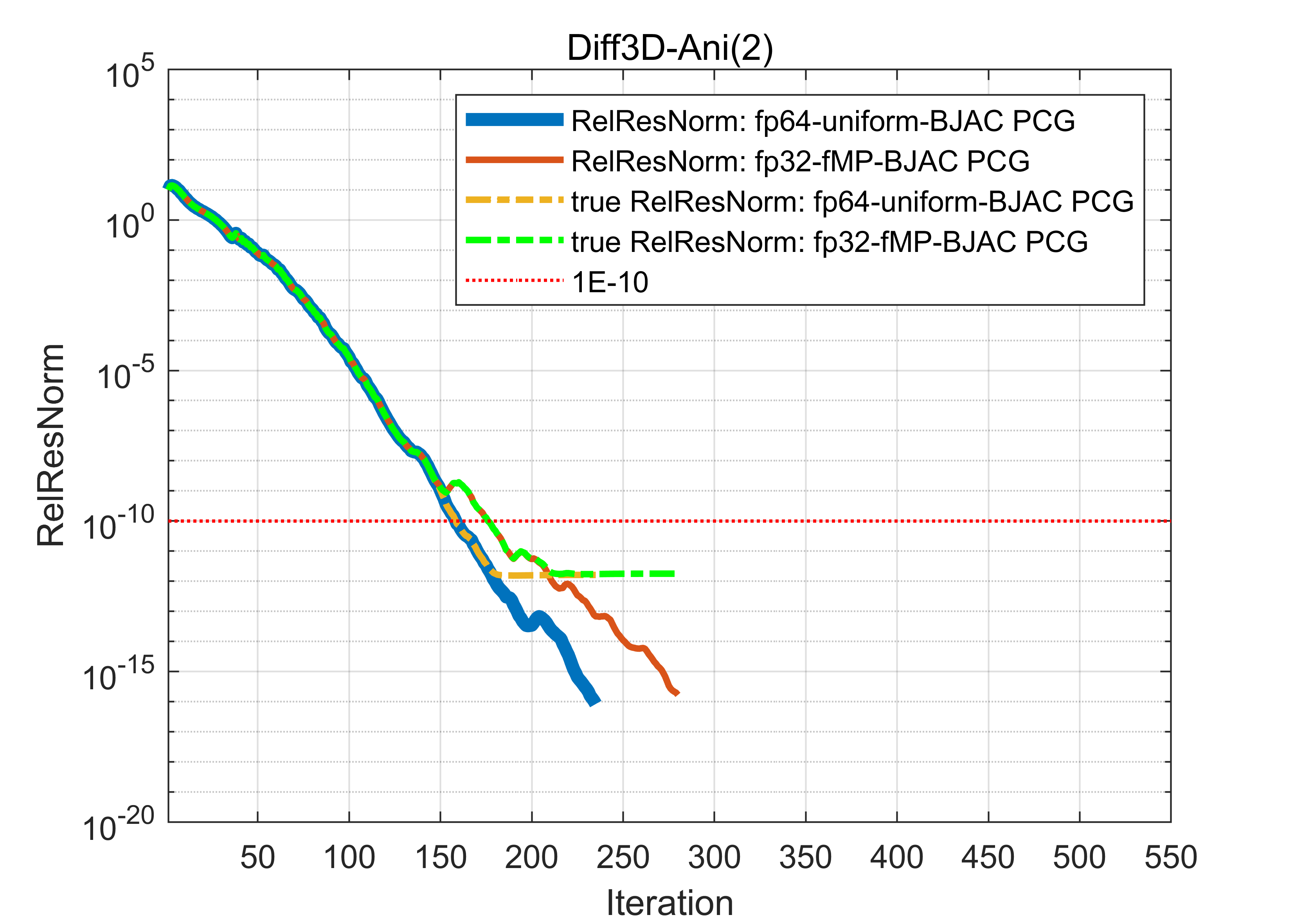}
    \includegraphics[width=0.47\linewidth]{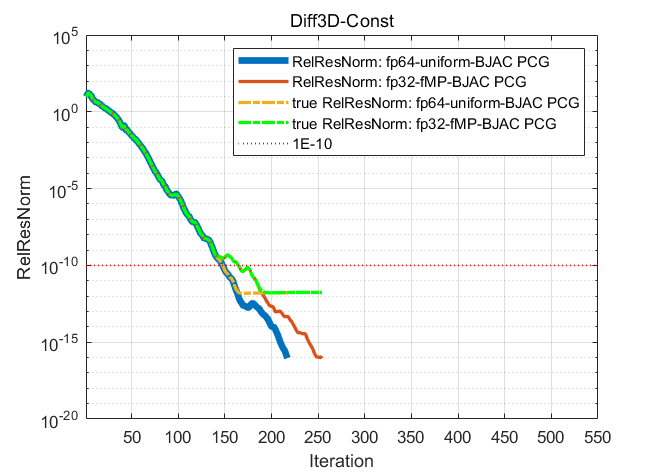} \\
    \caption{RelResNorm curves and true RelResNorm of the PCG algorithm for solving the Diff3D-Ani($s$) problem with $s = 1000, 100, 10, 4, 2, 1$.}
    
    \label{fig4.6}
\end{figure*}

Figure \ref{fig4.6} show the decline of relative residuals norm(RelResNorm) as well as the true RelResNorm for the uniform-BJAC and the fp32-fMP-BJAC with the six different $s$. Here, true RelResNorm uses residual computed via the residual definition formula $r=b-Ax$, instead of the residual deduced in PCG iterations.

From Figure \ref{fig4.6}, it can be found that for problem Diff3D-Ani(1000) , the two RelResNorm curves of the uniform-BJAC and the fp32-fMP-BJAC overlap completely, which indicates that the fp32-fMP-BJAC and the uniform-BJAC converge with the same convergence speed. For the problems Diff3D-Ani(100), Diff3D-Ani(10), Diff3D-Ani(4), Diff3D-Ani(2) and Diff3D-Ani(1), the two RelResNorm curves of the uniform-BJAC and fp32-fMP-BJAC are overlapped in the early iterations stage, but the bifurcations occur when the RelResNorms fall to $10^{-16}$, $10^{-11}$, $10^{-10}$, $10^{-9}$ and $10^{-10}$, respectively, and then the RelResNorm curves of fp32-fMP-BJAC suddenly rise and then fall again, finally reach to the convergence criterion. Since the convergence threshold is $10^{-10}$, for the problems Diff3D-Ani(1000), Diff3D-Ani(100), Diff3D-Ani(10), Diff3D-Ani(4), the number of converged iterations of the fMP-BJAC is the same as the uniform-BJAC, which are 327, 321, 205 and 179, respectively. It indicates that there is no delay phenomenon. However, for the Diff3D-Ani(2) and the Diff3D-Ani(1), i.e., Diff3D-Const, the number of iterations for these two preconditioners are different, 159 and 148 for uniform-BJAC while 177 and 165 for the fp32-fMP-BJAC, respectively, which show the convergence delay.

\begin{figure*}[!ht]
    \centering
    \includegraphics[width=0.47\linewidth]{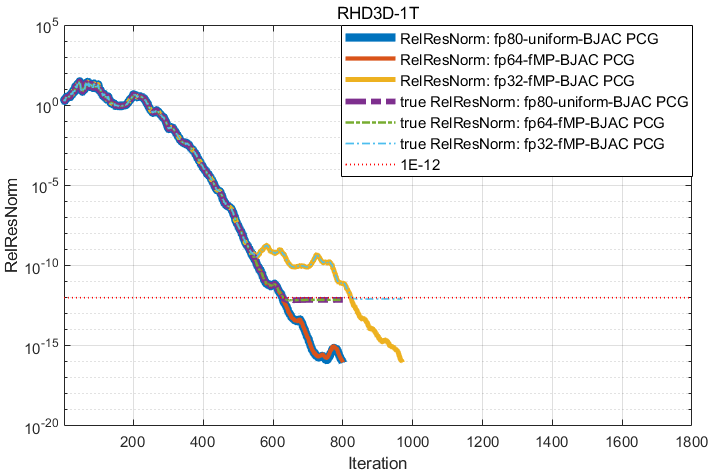} 
    \includegraphics[width=0.47\linewidth]{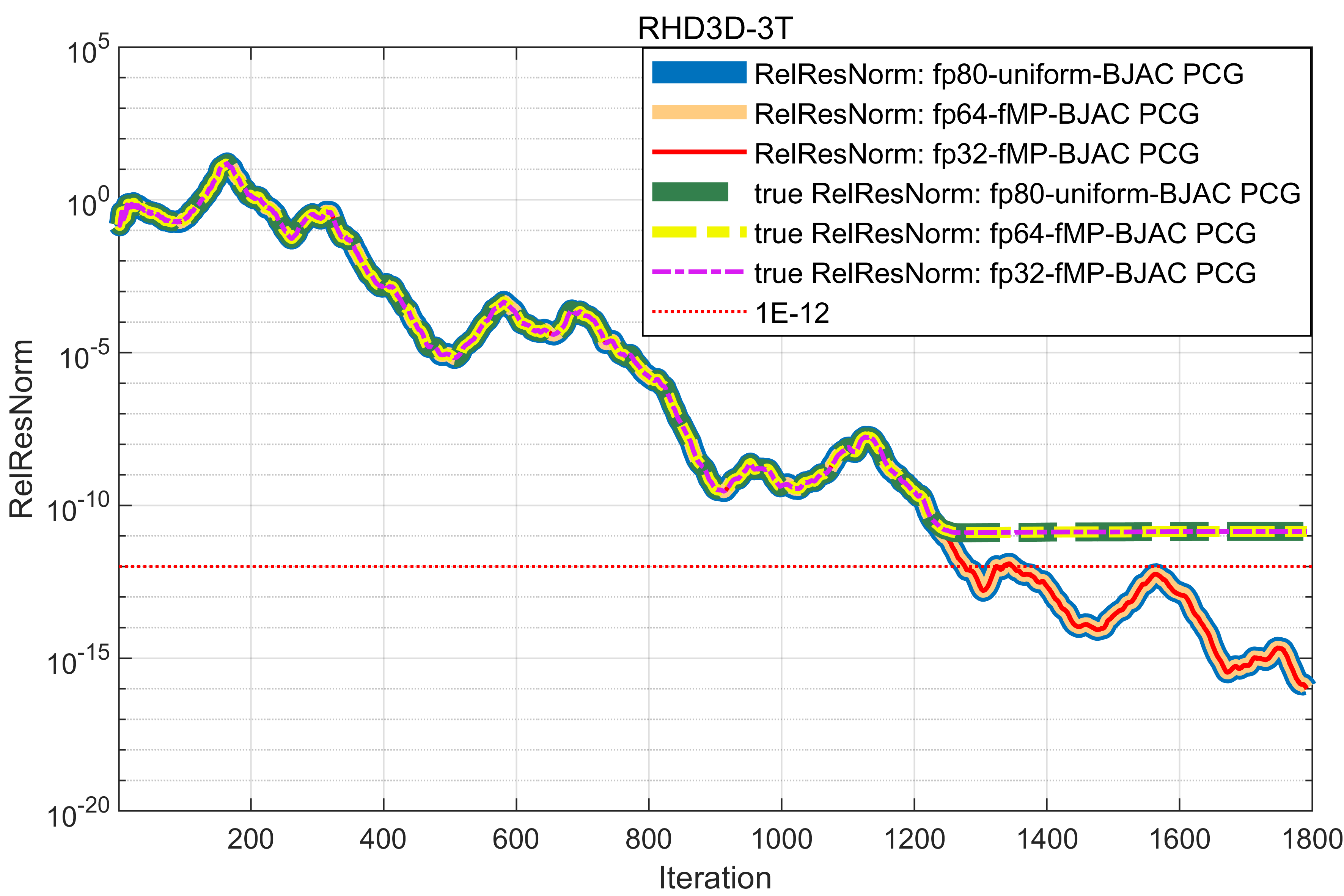} 
    \caption{RelResNorm and true RelResNorm curves of the PCG algorithm for solving the RHD3D-1T problem (left) and RHD3D-3T problem (right).}
    \label{fig4.7}
\end{figure*}

Similarly, for two practical problems, RHD3D-1T and RHD3D-3T, the multiscale strength distributions are shown in Table \ref{tab4.3}. It shows that the multiscale strength of RHD3D-1T is weaker than that of RHD3D-3T. The RelResNorm curves of the uniform-BJAC PCG, the fp64-fMP-BJAC PCG, and fp32-fMP-BJAC PCG are given in Figure \ref{fig4.7}. It can be found that for both RHD3D-1T and  RHD3D-3T, the RelResNorm curves of the uniform-BJAC and the fp64-fMP-BJAC overlap, which indicates that the fp64-fMP-BJAC does not lead to convergence delay. For RHD3D-1T, the RelResNorm curves of the uniform-BJAC and the fp32-fMP-BJAC bifurcate when RelResNorm fall to $10^{-9}$, and the number of iterations are 628 and 824, respectively. Since the convergence threshold is $10^{-12}$, that is, the fp32-fMP-BJAC leads to the convergence delay. For the RHD3D-3T problem, however, the RelResNorm curves of the uniform-BJAC and the fp32-fMP-BJAC overlap, indicating that the fp32-fMP-BJAC does not cause convergence delay.

It should be pointed out that, the true RelResNorm curves are also given in Figure \ref{fig4.6} and Figure \ref{fig4.7}, the results show that the reduction trend of the true RelResNorm curves are similar to the RelResNorm curves, but the difference is that the true RelResNorm curves for both preconditioners eventually flatten out to an almost identical convergence limit. It indicates that the fMP-BJAC and the uniform-BJAC have the same limiting accuracy.

\subsubsection{Diagonal Dominance Feature}\label{sec:4.5.2}

Diagonal Dominance is another common feature in many applications. A typical scene that leading to diagonal dominance is time-dependent problems, e.g. RHD3D-1T and RHD3D-3T problems considered in this paper, as shown in (\ref{eq4-2}). For time-dependent problems, 
the items of time discretization with a time step $dt$ will be add to the diagonal entries of matrix that lead to the diagonal dominance, and the strength of the diagonal dominance is closely related to the time step $dt$, the smaller $dt$, the stronger the diagonal dominance, and vice versa. Theoretical analyses \cite{Saad} show that the convergence of the Jacobi method is relevant to the strength of the diagonal dominance of the matrix. 

Here, taking RHD3D-1T as an example, we investigate the correlation between the convergence delay and diagonal dominance strength for mixed precision preconditioner. 
Specifically, for the row $i$, we change the diagonal dominance strength by adding a factor, denote as $P_{diag}$, of the sum of the absolute values of the off-diagonal entries to the diagonal entry:
\begin{math}
    a_{i,i}  \leftarrow a_{i,i} + P_{diag} \times \sum_{j=0,j\neq i}^{j=n} \vert a_{i,j} \vert
\end{math}
, where $a_{i,j}$ denotes the off-diagonal entry, $a_{i,i}$ denotes the diagonal entry. The larger the $P_{diag}$ is , the stronger the diagonal dominance is.

\begin{figure}[htb]
    \centering
    
    \includegraphics[width=1\linewidth]{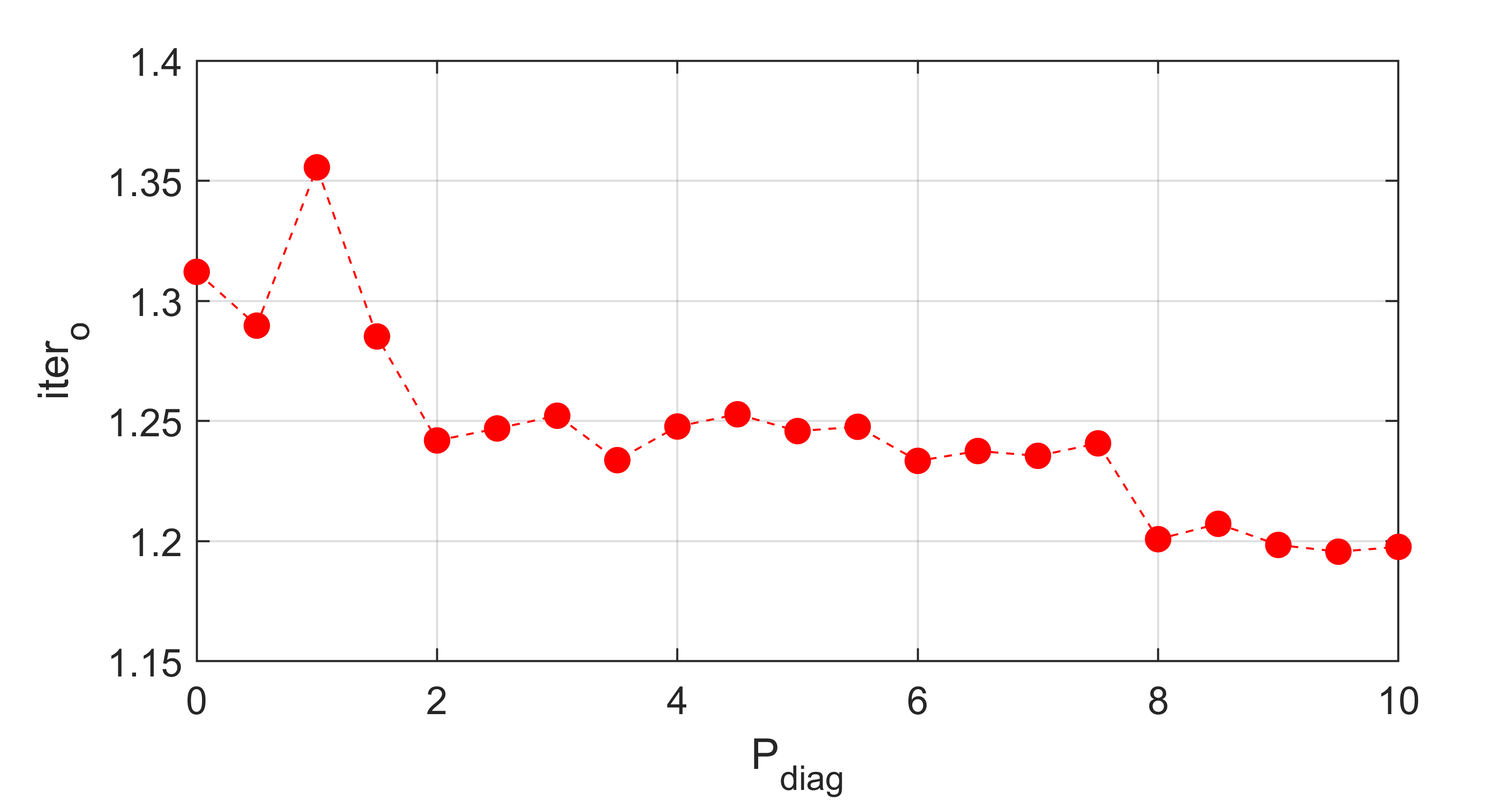}
    
    \caption{The changing of the overhead of additional iterations $iter_o$ with increasing the strength of diagonal dominance $P_{diag}$: RHD3D-1T problem.}
    \label{fig4.8}
\end{figure}

We observe the trend of changing in convergence delay as $P_{diag}$ increasing, the result is given in Figure \ref{fig4.8}, where $iter_{o} = iter_{m}/iter_{u}$ denotes the ratio of the iteration numbers of the fp32-fMP-BJAC to that of the uniform BJAC preconditioners, $iter_o$ reflects the overhead of additional iterations introduced by mixed precision preconditioner.

As shown in Figure \ref{fig4.8}, as the diagonal dominance strength increasing, $iter_{o}$ shows a fluctuating trend of declining. For example, the $iter_{o}$ is 1.36, 1.24, and 1.20 when $P_{diag}$ is 1, 2, and 8, respectively, suggesting that an increase in the diagonal dominance strength can improve the  convergence delay of the fMP-BJAC preconditioner.

\section{Conclusions}\label{sec:5}

%\subsection{Conclusions}\label{sec:5.1}

The mixed precision preconditioners fMP-BJAC and aMP-BJAC(hl) demonstrate significant potential for performance gains. These gains primarily depend on balancing the advantages of low precision computations and memory accesses for single iteration against the increase in additional iterations due to mixed precision. Experimental results on model problems and practical problems indicate that the fMP-BJAC and the aMP-BJAC(hl) can achieve speedup from $1.38 \times$ to $1.85 \times$ compared to the uniform high precision BJAC preconditioiner.

We find that, mixed precision preconditioners may cause convergence delays in some problems, i.e., additional iterations are required to achieve convergence, compared to the uniform high precision preconditioner. The adaptive mixed precision preconditioner aMP-BJAC(hl), which is based on the high-to-low precision order, can avoid convergence delays by tuning the adaptive threshold. The aMP-BJAC(lh) preconditioner which is based on the low-to-high precision order, however, can hardly avoid the convergence delays phenomenon. 

The numerical results show that the convergence delay behavior of the mixed precision preconditioner is related to the multiscale and diagonal-dominance 
features of matrix: the weaker the multiscale and the diagonal dominance, the more significant the convergence delay is. The theoretical analysis for the convergence behavior of the mixed precision preconditioners, as well as the optimal adaptive threshold determination for the aMP-BJAC preconditioner, are important issues which should be concerned in the future work.

\bmhead{Acknowledgements}

This research was funded by National Key Research and Development Program of China (No.2023YFB3001605), Science Challenge Project of China (No.TZ2024009), National Natural Science Foundation of China (No.62032023 and No.12301476).

\section*{Conflict of Interest}
On behalf of all authors, the corresponding author states that there is no confict of interest.

\bibliographystyle{bst/sn-mathphys-num}

\bibliography{CCFTHPC-bibliography}

\end{document}